\documentclass[a4paper,12pt]{amsart}
\usepackage{amsfonts}
\usepackage{amssymb}
\usepackage{ifthen}
\usepackage{graphicx}
\nonstopmode \numberwithin{equation}{section}
\usepackage[usenames]{color}
\newtheorem{thm}{Theorem}
\newtheorem{lem}{Lemma}
\newtheorem{cor}{Corollary}

\newtheorem{cl}{Claim}
\newtheorem{ca}{Case}
\newtheorem{sca}{Subcase}
\newtheorem{scl}{Subclaim}
\newtheorem{conj}[equation]{Conjecture}

\theoremstyle{definition}
\newtheorem{defn}{Definition}

\newtheorem{op}[equation]{Open Problem}
\newtheorem{ques}[equation]{Question}
\newtheorem{rem}{Remark}
\newtheorem{exam}[equation]{Example}

\newcounter {own}
\def\theown {\thesection       .\arabic{own}}

\newenvironment{pf}[1][]{%
 \vskip 3mm
 \noindent
 \ifthenelse{\equal{#1}{}}%
  {{\slshape Proof. }}%
  {{\slshape #1.} }%
 }%
{\qed\bigskip}

\newcounter{alphabet}
\newcounter{tmp}
\newenvironment{Thm}[1][]{\refstepcounter{alphabet}%
\medskip%
\noindent%
{\bf Theorem \Alph{alphabet}}%
\ifthenelse{\equal{#1}{}}{}{ (#1)}%
{\bf .} \itshape}{\vskip 8pt}

\makeatletter
\newcommand{\Ref}[1]{\@ifundefined{r@#1}{}{\setcounter{tmp}{\ref{#1}}\Alph{tmp}}}
\makeatother

\newcommand{\IR}{{\mathbb R}}

\newcommand{\IB}{{\mathbb B}}

\newcommand{\diam}{{\operatorname{diam}}}

\newcommand{\dist}{{\operatorname{dist}}}

\def\be{\begin{equation}}
\def\ee{\end{equation}}

\newcommand{\bee}{\begin{enumerate}}
\newcommand{\eee}{\end{enumerate}}

\newcommand{\blem}{\begin{lem}}
\newcommand{\elem}{\end{lem}}
\newcommand{\bthm}{\begin{thm}}
\newcommand{\ethm}{\end{thm}}
\newcommand{\bcor}{\begin{cor}}
\newcommand{\ecor}{\end{cor}}
\newcommand{\beg}{\begin{exam}}
\newcommand{\eeg}{\end{exam}}
\newcommand{\begs}{\begin{examples}}
\newcommand{\eegs}{\end{examples}}
\newcommand{\bdefe}{\begin{defn}}
\newcommand{\edefe}{\end{defn}}
\newcommand{\bprob}{\begin{prob}}
\newcommand{\eprob}{\end{prob}}
\newcommand{\bques}{\begin{ques}}
\newcommand{\eques}{\end{ques}}
\newcommand{\bei}{\begin{itemize}}
\newcommand{\eei}{\end{itemize}}
\newcommand{\bcon}{\begin{conj}}
\newcommand{\econ}{\end{conj}}
\newcommand{\bop}{\begin{op}}
\newcommand{\eop}{\end{op}}

\newcommand{\bca}{\begin{ca}}
\newcommand{\eca}{\end{ca}}
\newcommand{\bsca}{\begin{sca}}
\newcommand{\esca}{\end{sca}}

\newcommand{\bcl}{\begin{cl}}
\newcommand{\ecl}{\end{cl}}

\newcommand{\bscl}{\begin{scl}}
\newcommand{\escl}{\end{scl}}

\newcommand{\bcons}{\begin{conjs}}
\newcommand{\econs}{\end{conjs}}
\newcommand{\bprop}{\begin{propo}}
\newcommand{\eprop}{\end{propo}}
\newcommand{\br}{\begin{rem}}
\newcommand{\er}{\end{rem}}
\newcommand{\brs}{\begin{rems}}
\newcommand{\ers}{\end{rems}}
\newcommand{\bo}{\begin{obser}}
\newcommand{\eo}{\end{obser}}
\newcommand{\bos}{\begin{obsers}}
\newcommand{\eos}{\end{obsers}}
\newcommand{\bpf}{\begin{pf}}
\newcommand{\epf}{\end{pf}}
\newcommand{\ba}{\begin{array}}
\newcommand{\ea}{\end{array}}
\newcommand{\beq}{\begin{eqnarray}}
\newcommand{\beqq}{\begin{eqnarray*}}
\newcommand{\eeq}{\end{eqnarray}}
\newcommand{\eeqq}{\end{eqnarray*}}

\newcommand{\ds}{\displaystyle}

\newcounter{minutes}
\divide\time by 60
\newcounter{hours}
\multiply\time by 60 \addtocounter{minutes}{-\time}

\begin{document}

\bibliographystyle{amsplain}
\title{On quasisymmetry of quasiconformal mappings and its applications}

\thanks{
File:~\jobname .tex,
          printed: \number\year-\number\month-\number\day,
          \thehours.\ifnum\theminutes<10{0}\fi\theminutes}

\author{M. Huang}
\address{M. Huang, Department of Mathematics,
Hunan Normal University, Changsha,  Hunan 410081, People's Republic
of China} \email{mzhuang79@163.com}

\author{S. Ponnusamy $^\dagger $}
\address{S. Ponnusamy, Department of Mathematics,
Indian Institute of Technology Madras, Chennai-600 036, India.}
\email{samy@isichennai.res.in, samy@iitm.ac.in}

\author{A. Rasila}
\address{A. Rasila, Department of Mathematics and Systems Analysis, Aalto University,
FI-00076 Aalto, Finland} \email{antti.rasila@iki.fi}

\author{X. Wang ${}^{\dagger \dagger}$
}
\address{X. Wang, Department of Mathematics,
Hunan Normal University, Changsha,  Hunan 410081, People's Republic
of China} \email{xtwang@hunnu.edu.cn}

\keywords{Uniform domain, broad domain, John domain, $LLC_2$ domain, quasiconformal mapping, quasisymmetry.\\
$
^\dagger$ {\tt This author is currently at the
Indian Statistical Institute (ISI), Chennai Centre, SETS (Society
for Electronic Transactions and security), MGR Knowledge City, CIT
Campus, Taramani, Chennai 600 113, India.}\\
${}^{\dagger \dagger}$ {\tt Corresponding author}
}

\begin{abstract}
Suppose that $f:\, D\to D'$ is a quasiconformal mapping, where $D$ and $D'$ are domains in $\IR^n$, and that $D$ is a broad
domain. Then for every arcwise connected subset $A$ in $D$, the weak quasisymmetry of the restriction $f|_A:\, A\to f(A)$
implies its quasisymmetry, and as a consequence, we see that the answer to one of the open problems raised by
Heinonen from 1989 is affirmative under the additional condition that $A$ is arcwise connected.
As an application, we establish nine equivalent conditions for a bounded domain, which is quasiconformally
equivalent to a bounded and simply connected uniform domain, to be John. This result is a generalization of
the main result of Heinonen from \cite{JH}.
\end{abstract}


\maketitle \pagestyle{myheadings} \markboth{M. Huang, S. Ponnusamy, X. Wang, and
A. Rasila}{On quasisymmetry of quasiconformal mappings and its applications}

\section{Introduction and main results}\label{sec-1}

Quasisymmetric maps originate from the work of Beurling and Ahlfors \cite{BA}, who defined them
as the boundary values of quasiconformal self-maps of the upper half-plane onto the real line.
Since then this concept has proved to be very useful, and it has played a significant role,
e.g., in the work of Sullivan \cite{Su}. The general definition of quasisymmetry is due to
Tukia and V\"ais\"al\"a, who introduced the general class of quasisymmetric maps in \cite{TV},
and it has been studied by numerous authors thereafter, see for example \cite{BK, BoK, BoM, HK1, HK2, Ty, Vai0-1}.

In this paper, we study the quasisymmetry of quasiconformal mappings
in $\IR^n$ and certain applications of this property. Motivation for
this study arises from one of Heinonen's open problems together with
the main result, namely, Theorem $3.1$ of \cite{JH}. We now recall a result of
Heinonen, which is a generalization of a result of V\"ais\"al\"a \cite[Theorem
2.20]{Vai0}.

\begin{Thm}\label{Thm-1}$($\cite[Theorem 6.1]{JH}$)$ Suppose that $f:\, D\to D'$ is
a $K$-quasiconformal mapping, where $D$ and $D'$ are bounded domains in $\IR^n$, and that
$D$ is $\varphi$-broad. If $A\subset D$ is such that $f(A)$ is $b$-$LLC_2$
with respect to $\delta_{D'}$ in $D'$, then the
restriction $f|_A:\, A\to f(A)$ is weakly $H$-quasisymmetric in the metrics $\delta_{D}$ and
$\delta_{D'}$ with $H$ depending only on the
data
$$\omega=\left (n, K, b, \varphi, \frac{\delta_{D}(A)}{d_{D}(x_0)},
\frac{\delta_{D'}(f(A))}{d_{D'}(f(x_0))}\right ),
$$
where $x_0$ is some fixed point in $A$ and
$d_D(x_0)$ (resp. $\delta_{D}(A)$) denotes the distance from $x_0$ to the boundary $\partial D$ of $D$
(resp. $\delta_{D}$-diameter of $A$).
\end{Thm}

As a converse to Theorem \Ref{Thm-1}, Heinonen and  N\"akki \cite{JH-1}
further obtained the result below.

\begin{Thm}\label{Thm-2}$($\cite[Lemma 8.3]{JH-1}$)$ Suppose that $f:\, D\to D'$ is
a $K$-quasiconformal mapping, where $D$ and $D'$ are proper domains in $\IR^n$, and that $D'$ is $\varphi$-broad.
If $A\subset D$ is arcwise connected and
$f|_A:\,A\to f(A)$ is weakly $H$-quasisymmetric in the metrics
$\delta_{D'}$ and $\delta_D$, then
$A$ is $b$-$LLC_2$ with respect to
$\delta_{D}$ in $D$, where $b$ depends only on the
data
$$\varpi=(n, K, \varphi, H).
$$
\end{Thm}

In \cite{Vai0}, V\"{a}is\"{a}l\"{a} proved that every weak quasisymmetry
$f:\, X\to Y$ is quasisymmetric provided that both $X$ and $Y$ are $HTB$ metric
spaces and that $X$ is arcwise connected (see \cite{TV} for the definition of $HTB$ spaces).
Heinonen pointed out in \cite{JH} that this amenable $HTB$-criterion is not automatically
satisfied as there are domains which are $LLC_2$ with respect to $\delta_D$, but which are
not $HTB$. Hence, Heinonen asked whether the condition ``weakly" in Theorem \Ref{Thm-1}
is redundant or not (see the paragraph next to the statement of \cite[Theorem $6.5$]{JH}).
In this paper, we first study this problem. Our result is the following.

\begin{thm}\label{thm1}
Suppose that $f:\,
D\to D'$ is a $K$-quasiconformal mapping, where  $D$ and $D'$ are  domains in $\IR^n$,
and that $D$ is $\varphi$-broad. For an arcwise connected set
$A$ in $D$, if the restriction $f|_A:\, A\to f(A)$ is
weakly $H$-quasisymmetric in the metrics
$\delta_{D}$ and $\delta_{D'}$, then
$f|_A:\, A\to f(A)$ is $\eta$-quasisymmetric in the metrics
$\delta_{D}$ and $\delta_{D'}$ with $\eta$ depending only on the
data
$$\mu=(n,K,H,\varphi).
$$
\end{thm}

The next result easily follows from Theorems \ref{thm1},
\Ref{Thm-1} and \Ref{Thm-2} together with Remark \ref{xt-11} given in Section \ref{sec-2}.





 \bcor\label{Sat-4} Suppose that $D$ and $D'$ are bounded domains in $\IR^n$, that  $f:\, D\to D'$ is a
$K$-quasiconformal mapping, where $D$ is a $\varphi$-broad domain, and that $A$
is an arcwise connected subset of $D$. Then the following statements
are equivalent.

\begin{enumerate}

\item[{\rm (1)}]\label{5-4}
$f|_A:\, A\to f(A)$ is $\eta$-quasisymmetric in the
metrics $\delta_{D}$ and $\delta_{D'}$;

\item[{\rm (2)}]\label{5-1}
$f(A)$ is $b$-$LLC_2$ with respect to $\delta_{D'}$ in $D'$;

\item[{\rm (3)}]\label{5-3}
$f|_A:\, A\to f(A)$ is weakly $H$-quasisymmetric in the
metrics $\delta_{D}$ and $\delta_{D'}$,
\end{enumerate}
where $b$, $H$
and $\eta$ depend on each other and the data
$$u=\left (n,K,\varphi, \frac{\delta_{D}(A)}{d_{D}(x_0)},
\frac{\delta_{D'}(f(A))}{d_{D'}(f(x_0))}\right ),
$$
and $x_0$ is a fixed point in $A$.
 \ecor

\br
\begin{enumerate}

\item[{\rm (1)}] The equivalence of $(1)$ and $(2)$ in Corollary \ref{Sat-4} shows that the answer to Heinonen's
problem mentioned as above is affirmative when the set $A$ is
arcwise connected;

\item[{\rm (2)}]The equivalence of $(2)$ and $(3)$ in Corollary \ref{Sat-4} shows that the converse of Theorem \Ref{Thm-1} is also true when the set $A$ is
arcwise connected.
\end{enumerate}
\er

In \cite{JH}, Heinonen studied the quasiconformal mappings of the
unit ball $\mathbb{B}$ in $\IR^n$ onto John domains $D$ in $\IR^n$.
The main results of \cite{JH} provide nine equivalent conditions for a bounded domain $D$,
which is quasiconformally equivalent to $\mathbb{B}$, to be John, see \cite[Theorem
3.1]{JH}. In
addition, Heinonen specially pointed out that the requirement that ``$D$
is quasiconformally equivalent to $\mathbb{B}$" in \cite[Theorem
3.1]{JH} cannot be replaced, e.g., by the requirement that ``$D$ is homeomorphic to
$\mathbb{B}$" or ``$D$ is a Jordan domain". In this paper,
we shall further refine this result. Based on Theorem
\ref{thm1}, we shall actually prove Theorem \ref{thm2} below, which
shows that the ball ``$\IB$" in the requirement that ``$D$ being
quasiconformally equivalent to $\mathbb{B}$" in \cite[Theorem
3.1]{JH} can be replaced by ``a bounded and simply connected uniform domain". Note that
every ball in $\IR^n$ is a bounded and simply connected uniform domain. To state the result, we first recall some notations.

Suppose $D$ is a bounded and simply connected $c$-uniform
domain. For $x\in D$, we use $\Phi(x)$ to denote the set of all
components $I(x)$ in the intersection $\mathbb{B}(x, 8cd_D(x))\cap
\partial D$ such that $\diam(I(x))\geq d_D(x)$, where ``$\diam$" means ``diameter".
Obviously, $\Phi(x)\not=\emptyset$ for each $x\in D$, but it is possible that
$\Phi(x)$ contains only one element.

\begin{thm} \label{thm2}
Suppose that $D$ and $D'$ are bounded domains in $\IR^n$, that $f:\,
D\to D'$ is a $K$-quasiconformal mapping, that $f:\,\overline{D}\to
\overline{D'}$ is continuous, and $x_0\in D$. If $D$ is a simply
connected $c$-uniform domain, then the following statements are
equivalent.
\begin{enumerate}
\item\label{thm2-1} $D'$ is a $b$-John domain with center $f(x_0)$;

\item\label{thm2-2} $D'$ is $\varphi$-broad;

\item\label{thm2-3} $f:\, (D,\delta_D)\to (D', \delta_{D'})$ is
$\eta$-quasisymmetric;

\item\label{thm2-4} For $x\in D$ and each
$I(x)\in \Phi(x)$, $\diam\big(f(I(x))\big)\leq b_1d_{D'}(f(x))$;

\item\label{thm2-5} For $x, w \in D$, if $|x-w|\leq 8cd_D(x)$, then
$\delta_{D'}(f(x),f(w))\leq b_2d_{D'}(f(x))$;

\item\label{thm2-6} For $x, w \in D$, if $|x-w|\leq 8cd_D(x)$ and $d_D(w)\leq
2cd_D(x)$, then
$$a_f(w)\leq b_3a_f(x)\Big(\frac{d_D(x)}{d_D(w)}\Big)^{1-\alpha};
$$

\item\label{thm2-7} For all components
$P\subset Q$, where $P\in \Phi(x)$ and $Q\in \Phi(w)$,
\[
\ds \frac{\diam(f(P))}{\diam(f(Q))}\leq
b_4\Big(\frac{\diam(P)}{\diam(Q)}\Big)^{\alpha}
\]
for $x$, $w\in D$ (here the case $x=w$ is included);

\item\label{thm2-8} $D'$ is $b_5$-$LLC_2$;

\item\label{thm2-9} $D'$ is $b_6$-$LLC_2$ with respect to $\delta_{D'}$;

\item\label{thm2-10} $f:\, (D, \delta_{D})\to (D', \delta_{D'})$ is weakly $H$-quasisymmetric,
\end{enumerate}
where the constants $b$, $b_1$, $b_2$, $b_3$, $b_4$, $b_5$, $b_6$,
$\alpha$, $H$ and the functions $\varphi$, $\eta$ depend only on each
other and the data
$$v=\Big(n, K, c,\frac{\diam(D)}{d_D(x_0)},
\frac{\diam(D')}{d_{D'}(f(x_0))} \Big).
$$
\end{thm}

This paper is organized as follows. In Section \ref{sec-2},
we shall introduce necessary notations, and recall some
preliminary results. We shall prove Theorem \ref{thm1} and Theorem \ref{thm2} in Section \ref{sec-3} and \ref{sec-4}, respectively. The proofs are mainly based on the properties of the conformal modulus of a curve family.

\section{Preliminaries}\label{sec-2}

\subsection{Notation}
Throughout the paper, we always assume that $D$ and $D'$ are domains in $\IR^n$, $n\geq 2$, and
$f:\, D\to D'$ includes the assumption that $f$ is a homeomorphism  from $D$ onto $D'$.
Also we use $\mathbb{B}(x_0,r)$ to denote the open ball $\{x\in
\IR^n:\,|x-x_0|<r\}$ centered at $x_0$ with radius $r>0$. Similarly,
for the closed balls and spheres, we use the  notations
$\overline{\mathbb{B}}(x_0,r)$ and $\mathbb{S}(x_0,r)$,
respectively. In particular, we use $\IB$ to denote the unit ball
$\mathbb{B}(0,1)$ and $\mathbb{S}$ its boundary.

For convenience, in what follows, we always assume that $x, y, z,
\ldots$ are points in $D$ and the primes $x',
y', z', \ldots$ denote the images of $x, y, z, \ldots$ in $D'$ under $f$,
respectively. Also we assume that $\alpha, \beta, \gamma, \ldots$
are curves in $D$ and the primes $\alpha', \beta', \gamma', \ldots$ denote
the images of $\alpha, \beta, \gamma, \ldots$ in $D'$ under $f$,
respectively. For a set $A$ in $D$, $A'$ denotes the image of $A$ in $D'$ under $f$.

\subsection{John domains and uniform domains}
John \cite{Jo}, Martio and Sarvas \cite{MS}  were the first who
studied John domains and uniform domains, respectively. There are
many alternative characterizations for uniform and John domains, see
\cite{Bro, FW,  Geo, Kil, Martio-80, Vai4, Vai5, Vai6, Vai7, Vai8}.
The importance of these two classes of domains in the function theory is
well documented, see \cite{FW, Kil, RJ, Vai2}. Moreover, John
domains and uniform domains in $\mathbb{R}^n$ have numerous
geometric and function theoretic features, which are useful in many areas of modern
mathematical analysis, see \cite{Alv, Bea, Bro, Geo, Has, Ml,
M2, Jo80, Jo81, Yli, Vai2, Vai7}. From the various equivalent
characterizations, we adopt the following definitions.

\bdefe \label{def1}
A domain $D$ in $\IR^n$  is said to be  {\it $c$-uniform}  if there
exists a constant $c$ with the property that each pair of points
$z_{1},z_{2}$ in $D$ can be joined by a rectifiable arc $\gamma$ in
$ D$ satisfying (cf. \cite{MS, Vai4})
\begin{enumerate}
\item\label{eq-1}
$\ds\min_{j=1,2}\ell (\gamma [z_j, z])\leq c\, d_D(z)$ for all $z\in \gamma$, and

\item\label{eq-2}
$\ell(\gamma)\leq c\,|z_{1}-z_{2}|$,
\end{enumerate}
where $\ell(\gamma)$ denotes the arc length of $\gamma$,
$\gamma[z_{j},z]$ the part of $\gamma$ between $z_{j}$ and $z$.
Also we say that $\gamma$ is a {\it double $c$-cone arc}.

A domain $D$ in $\IR^n$  is said to be a {\it $c$-John domain} if it
satisfies the condition  \eqref{eq-1} in Definition \ref{def1}, but
not necessarily \eqref{eq-2}. In this case, $\gamma$ is called a {\it $c$-cone arc}.
\edefe

\bdefe\label{wed-1}
A domain $D$ in $\IR^n$  is said to have the
{\it $c$-carrot property} with center $x_0\in \overline{D}$ if there
exists a constant $c$ with the property that for each point $z_{1}$
in $D$, $z_1$ and $x_0$ can be joined by a rectifiable arc $\gamma$
in $ D$ satisfying (cf. \cite{RJ, Vai0})
$$\ell (\gamma [z_1, z])\leq c\, d_D(z)
$$
for all $z\in \gamma$. Also we say that $\gamma$ is a {\it
$c$-carrot arc}.

A domain $D$ in $\overline{\IR}^n$  is said to be a {\it $c$-John
domain with center $x_0$ in $\overline{D}$}  if it has the
$c$-carrot property with center $x_0\in\overline{D}$.
\edefe

Definition \ref{wed-1} is often referred to as the ``arc length"
definition for the carrot property (resp. John domains). When the word ``arc length" in
Definition \ref{wed-1} is replaced by ``diameter", then it is called
the ``diameter" definition for the carrot property (resp. John domains). The following
result reveals the close relationship between these two definitions.

\begin{Thm}\label{ThmF-1}$($\cite{RJ}$)$
The ``arc length" definition for the carrot property or John domains is
quantitatively equivalent to the ``diameter" one with the same center.
\end{Thm}

Also the following result concerning the equivalence of the definitions for John
domains in Definitions \ref{def1} and \ref{wed-1} is due to V\"ais\"al\"a.

\begin{Thm}\label{sat-1} $($\cite[Lemma 2.4]{Vai0}$)$
The definitions for John domains in Definitions \ref{def1} and
\ref{wed-1} are quantitatively equivalent for bounded domains.
\end{Thm}

\subsection{Quasihyperbolic metric, solid arcs and linearly locally connected sets}
Let $\gamma$ be a rectifiable arc or path in $D$. Then the {\it
quasihyperbolic length} of $\gamma$ is defined to be the number
$\ell_{k_D}(\gamma)$ given by (cf. \cite{GP})
$$\ell_{k_D}(\gamma)=\int_{\gamma}\frac{|dz|}{d_D(z)}.
$$

For $z_1$, $z_2$ in $D$, the {\it quasihyperbolic distance}
$k_D(z_1,z_2)$ between $z_1$ and $z_2$ is defined in the usual way:
$$k_D(z_1,z_2)=\inf\ell_{k_D}(\gamma),
$$
where the infimum is taken over all rectifiable arcs $\gamma$
joining $z_1$ to $z_2$ in $D$. An arc $\gamma$ from $z_1$ to $z_2$
is called a {\it quasihyperbolic geodesic} if
$\ell_{k_D}(\gamma)=k_D(z_1,z_2)$. Each subarc of a quasihyperbolic
geodesic is obviously a quasihyperbolic geodesic. It is known that a
quasihyperbolic geodesic between two points in $D$ always exists
(cf. \cite[Lemma 1]{Geo}). Moreover, for $z_1$, $z_2$ in $D$, we
have (cf. \cite{Vai4, Vu-book-88})
\beq\label{eq-2-9'}
k_{D}(z_1, z_2) &\geq & \inf_{\gamma}\,\log\left
(1+\frac{\ell(\gamma)}{\min\{d_D(z_1), d_D(z_2)\}}\right )\\
\nonumber &\geq & \log\left (1+\frac{|z_1-z_2|}{\min\{d_D(z_1),
d_D(z_2)\}}\right ) \\
\nonumber &\geq &  \Big|\log \frac{d_D(z_2)}{d_D(z_1)}\Big|,
\eeq
where $\gamma$ denote rectifiable curves joining $z_1$ and
$z_2$ in $D$. In particular, it follows that for every quasigeodesic
$\gamma$ in $D$ joining $z_1$ to $z_2$,
\beq\label{eq-2-9}
k_{D}(z_1, z_2)\geq \log\left (1+\frac{\ell(\gamma)}{\min\{d_D(z_1),
d_D(z_2)\}}\right).
\eeq

 \noindent Furthermore, if $|z_1-z_2|\le d_D(z_1)$, then we have
(cf. \cite{Vai6-0, Vu})
\begin{equation} \label{upperbdk}
k_D(z_1,z_2)\le \log\Big( 1+ \frac{
|z_1-z_2|}{d_D(z_1)-|z_1-z_2|}\Big).
\end{equation}

The following characterization of uniform domains by the
quasihyperbolic metric is useful for our discussions.

\begin{Thm}\label{ThmF'} $($\cite[2.50 (2)]{Vu}$)$\quad
A domain $D\subset \IR^n$ is $c$-uniform if and only if there is a
constant $\mu_1$ such that for all $x$, $y\in D$,
$$k_{D}(x, y)\leq \mu_1\log \left(1+\frac{|x-y|}{\min\{d_{D}(x), d_{D}(y)\}}\right ),
$$
where $\mu_1=\mu_1(c)$ which means that $\mu_1$ is a constant depending only on $c$.
\end{Thm}

This form  of the definition of uniform domains is due to
Gehring and Osgood \cite{Geo}. As a matter of fact, in \cite[Theorem
1]{Geo}, there was an additive constant in the inequality of Theorem
\Ref{ThmF'}, but it was shown by Vuorinen in \cite[2.50(2)]{Vu}
that the additive constant can be chosen to be zero.

Next, we recall a relationship between the quasihyperbolic distance of
points in $D$ and the one of their images in $D'$ under a
quasiconformal mapping.

\begin{Thm} \label{ThmF}$($\cite[Theorem 3]{Geo}$)$
Suppose $f:\,D\to D'$ is a $K$-quasiconformal mapping. Then for $z_1$, $z_2\in
D$,
$$k_{D'}(z'_1,z'_2)\leq \mu_2\max\big\{k_{D}(z_1,z_2), (k_{D}(z_1,z_2))^{\frac{1}{\mu_2}}\big\},
$$
where $\mu_2=\mu_2(n, K)\geq 1$.
\end{Thm}






\bdefe \label{def2}
Suppose that $A\subset D$ and $b\geq 1$ is a
constant.  We say that $A$ is $b$-$LLC$$_2$ (resp. {\it $b$-$LLC$$_2$
with respect to $\delta_{D}$}) in $D$ if for all $x\in A$ and $r>0$,
the points in $A \backslash \overline{\mathbb{B}}(x, br)$ $($resp.
$A \backslash \overline{\mathbb{B}}_{\delta_{D}}(x, br))$ can be
joined in $D\backslash \overline{\mathbb{B}}(x, r)$ $($resp.
$D\backslash \overline{\mathbb{B}}_{\delta_{D}}(x, r))$, where
$$\overline{\mathbb{B}}_{\delta_{D}}(x, r)=\{z\in D:\, \delta_{D}(z, x)< r\}.
$$

If $A=D$, then we say that $D$ is $b$-\it{LLC}$_2$ $($resp.
$b$-\it{LLC}$_2$ with respect to $\delta_{D})$.
\edefe

\subsection{Moduli of families of curves}

Suppose that $G$ is a domain in $\IR^n$, that $E$ and $F$ are two
disjoint continua in $G$, and that ${\rm Mod}\,(E,F;G)$ denotes the
usual conformal modulus of the family of all curves joining $E$ and
$F$ in $G$. For a family of  curves $\Gamma$ in $G$, we always
 use ${\rm Mod}\,(\Gamma)$ to denote the conformal modulus of $\Gamma$
\cite{Vai}. The following related results are useful for us. The
first result is from \cite[p.~397]{HR} or the combination of \cite[\S
11.9]{Vai}, \cite[Lemmas 2.39 and 2.44]{Vu} and
\cite[\S 7]{Vu-book-88}.

\begin{Thm}\label{ThmG} Suppose $n \ge 2$. Then there
exist decreasing homeomorphisms  $\phi_n$, $\psi_n:\, (0, \infty) \to
(0, \infty)$ such that
$$\phi_n(t)\leq {\rm Mod}\,(E,F;\IR^n)\leq \psi_n(t),
$$
where $t=\frac{\dist(E,F)}{\min\{\diam E, \diam F\}}$ and ``dist"
means ``distance".
\end{Thm}

\begin{Thm}  \label{lem0-0}
{\rm (\cite[Theorem 4.15]{Gem}  and \cite[p.~397]{HR})} Suppose
that $G\subset \IR^n$ is a $c$-uniform domain. Then
$${\rm Mod}\,(E,F;\IR^n)\leq \mu_4{\rm Mod}\,(E,F;G)
$$
for every pair of disjoint continua $E$ and $F$ in $G$, where
$\mu_4=\mu_4(n,c)$.
\end{Thm}

\begin{Thm}\label{lemm1--0''} $($\cite[Theorem 7.1]{Vai} and \cite[Lemma 2.9]{Vai11}$)$
$(1)$ There is a decreasing homeomorphism $\varrho_n:\, (0,\infty)\to
(0,\infty)$ with the following property: if $\Gamma$ is a family of
paths, each of which meets a set $E$ in $\IR^n$ and has length at
least $\lambda$, then
$${\rm Mod}(\Gamma)\leq \varrho_n\Big(\frac{\lambda}{\diam(E)}\Big).
$$

$(2)$ Suppose that a family of the curves $\Gamma$ lie in a Borel
set $E\subset\overline{\IR}^n$ and that $\ell(\gamma)\geq r>0$ for
every locally rectifiable $\gamma\in \Gamma$. Then
$${\rm Mod}_p(\Gamma)\leq \frac{m(E)}{r^p},
$$
where $2\leq p\leq n$ and $m(E)$ denotes the volume of $E$.
\end{Thm}

\begin{Thm}\label{lemm1--xx} $($\cite[Section 7.5]{Vai}$)$ For $x\in \IR^n$ and
$0<a<b<\infty$, let $A$ denote the spherical ring
$\mathbb{B}(x,b)\backslash \mathbb{\overline{B}}(x,a)$,
$E=\mathbb{S}(x,a)$, $F=\mathbb{S}(x,b)$ and
$\Gamma_{A}=\Gamma(E,F; A)$ the family of the curves in $A$
connecting $E$ and $F$. Then
$${\rm Mod}(\Gamma_{A})=\omega_{n-1}\Big(\log \frac{b}{a}\Big)^{1-n},
$$
where $\omega_{n-1}$ denotes the $(n-1)$-dimensional surface area of $\mathbb{S}$.
\end{Thm}

\subsection{Internal metric, broad domains and quasisymmetric
mappings}

For $x$, $y$ in $D$, the internal metric $\delta_{D}$
in $D$ is defined by
$$\delta_D(x,y)=\inf \{\diam(\alpha):\, \alpha\subset D\;
\mbox{is a rectifiable arc joining}\; x\; \mbox{and}\; y \}.
$$

\bdefe \label{def3}
Let $\varphi:\, (0, \infty)\to  (0, \infty)$ be a
decreasing homeomorphism. We say that $D$ is {\it $\varphi$-broad}
if for each $t>0$ and each pair $(C_0, C_1)$ of continua in $D$
 with $C_0\cap C_1=\emptyset$, the condition
$\delta_{D}(C_0, C_1)\leq t\min\{\diam(C_0),\diam(C_1)\}$ implies
$${\rm Mod}\,(C_0, C_1; D)\geq \varphi(t),
$$
where $\delta_{D}(C_0, C_1)$ denotes the $\delta_{D}$-distance
between $C_0$ and $C_1$.
\edefe

Broad domains were introduced in \cite{Vai0}. It was later proved
that a simply connected planar domain is broad if and only if it is
John  \cite[Section 8]{RJ}. Further, Gehring and Martio proved that
each uniform domain in $\IR^n$ is
broad, see \cite[Lemma 2.6]{Gem}.

It is important to recall that the notion of broad domains also goes under the
term L\"{o}ewner space.
The notion of a Loewner space was introduced by Heinonen and Koskela \cite{HK2} in their
study of quasiconformal mappings of metric spaces; Heinonen's recent monograph \cite{Hei} renders an enlightening account of these ideas. See \cite{BK, BHK, Her, Ty} etc for more related discussions.

\bdefe \label{def2'}
Let $(X, d)$ and $(X', d')$ be two metric spaces, and let $\eta:\, [0, \infty)\to [0,
\infty)$ be a homeomorphism. An embedding $f:\, X\to X'$ is {\it
$\eta$-quasisymmetric}, or briefly $\eta$-QS, in the metrics $d$ and $d'$ if $d(a,x)\leq td(a,y)$ implies
$$d'(a',x')\leq \eta(t)d'(a',y')
$$
for all $a$, $x$, $y\in X$, and if  there is a constant $\nu\geq 1$ such
that $d(a,x)\leq d(a,y)$ implies
$$d'(a',x')\leq \nu d'(a',y'),
$$
then $f$ is said to be {\it weakly $\nu$-quasisymmetric}, or briefly weakly $\nu$-QS,
in the metrics $d$ and $d'$.
\edefe

Obviously, ``quasisymmetry" implies ``weak quasisymmetry".

\br\label{xt-11} It follows from  \cite{BHK}
that $f^{-1}$ is $\eta_1$-quasisymmetric if $f$ is
$\eta$-quasisymmetric, where $\eta_1(t)=\frac{1}{\eta^{-1}(1/t)}$.
It also follows from \cite[Lemma 3.9]{RJ} (resp. \cite[Lemma 3.5]{RJ}) that QS mappings preserve broad domains (resp. John domains).
\er

\subsection{The function $a_f$}

\bdefe \label{def2"}
Suppose $f:\,D\to D'$ is a $K$-quasiconformal
mapping. For $x\in D$, we write
$$B_x=\mathbb{B}\Big(x,\frac{1}{2}d_D(x)\Big)
$$
and set
$$a_f(x)=\exp\Big({\frac{1}{nm(B_x)}\int_{B_x}\log J_f\,dm}\Big),
$$
where $J_f$ denotes the Jacobian of $f$ and
$m(B_x)$ stands for the volume of the ball $B_x$.
\edefe

We recall the following result concerning the function $a_f$.

\begin{Thm}\label{mzxl-3-1} $($\cite[Lemma 2.11]{JH}, see also \cite[Theorem 1.8]{AG}$)$
Suppose $f:\,D\to D'$ is a $K$-quasiconformal mapping. Then there is a
constant $\mu_5$ such that
$$\frac{1}{\mu_5}\frac{d_{D'}(x')}{d_D(x)}\leq a_f(x)\leq \mu_5\frac{d_{D'}(x')}{d_D(x)}
$$
for all $x\in D$, where $\mu_5=\mu_5(n,K)$.
\end{Thm}

\section{ Linear local connectedness and quasisymmetry }\label{sec-3}

The aim of this section is to give a proof of Theorem \ref{thm1}.
Before the proof of Theorem \ref{thm1}, we establish two auxiliary lemmas.

\begin{lem}\label{Lemc1--0-00} Suppose that $f:\, D\to D'$ is a
$K$-quasiconformal mapping, that $D$ is $\varphi$-broad, and that
$A\subset D$ is arcwise connected such that the restriction $f|_A:\, A\to A'$ is
weakly $H$-QS in the metrics $\delta_{D}$ and $\delta_{D'}$. For
$z_1, z_2, z_3\in A$, if $\delta_D(z_1,z_3)\leq
c\delta_D(z_1,z_2)$, where $c$ is a constant, then
$$\delta_{D'}(z'_1,z'_3)\leq \mu_7(c)\delta_{D'}(z'_1,z'_2),
$$
where $\mu_7(c)=H\mu_6^{1+\log_3c}$,
\[
\mu_6=\max\bigg\{8HK,8H\exp\bigg({2\Big(\frac{2K\omega_{n-1}}
{\varphi(45/4)}\Big)^{\frac{2}{n-1}}}\bigg)
+8H, 2H(H+1)\bigg(\varrho_n^{-1}\Big(\frac{\varphi(8)}{K}\Big)+1\bigg)\bigg\},
\]
and $\varrho_n$ is the same as in Theorem {\rm \Ref{lemm1--0''}}.
\end{lem}

\bpf If $\delta_D(z_1,z_2)\geq \delta_D(z_1,z_3)$, then the assumption ``$f|_A$ being
weakly $H$-QS in the metrics $\delta_{D}$ and $\delta_{D'}$" implies
\beq\label{lawu-1} \delta_{D'}(z'_1,z'_3)\leq H\delta_{D'}(z'_1,z'_2).
\eeq

For the other case, that is, $\delta_D(z_1,z_2)< \delta_D(z_1,z_3),$
we take $\beta$ to be an arc joining $z_1$ and $z_3$ in $A$.
We partition $\beta$ with the aid of a finite sequence $\{v_i\}_{i=1}^{s+2}$ of points in $\beta$ as follows.

If $\delta_D(z_1,z_3)<3\delta_D(z_1,z_2)$, then we let $s=1$ and $x_s=z_1$.

\begin{figure}[htbp]
\begin{center}
\input{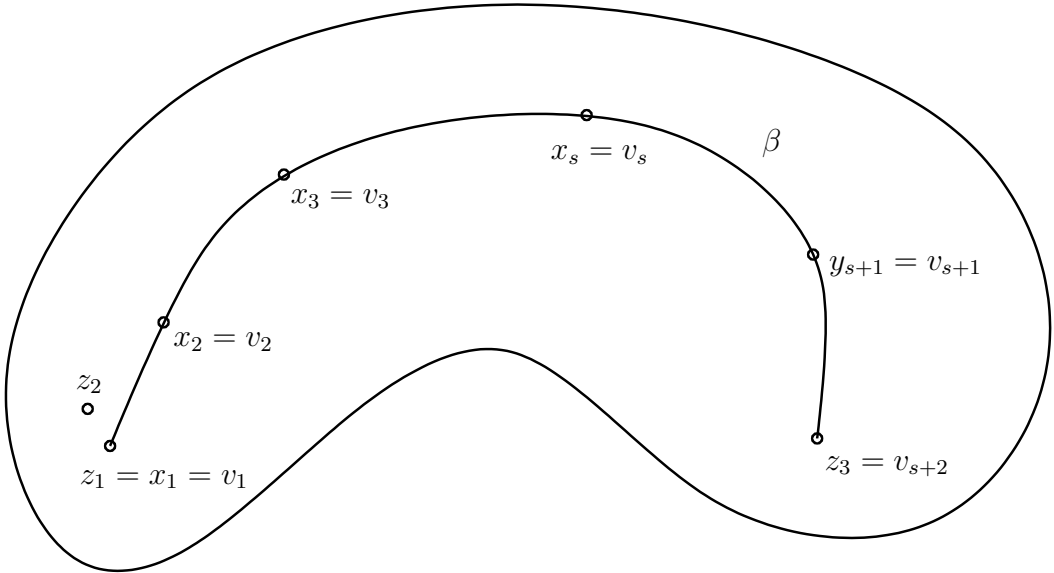}
\caption{The arc $\beta$ in $A$ and its partition}

\label{hconj-fig01}
\end{center}
\end{figure}

If $\delta_D(z_1,z_3)\geq3\delta_D(z_1,z_2)$, let $s$ be the number which satisfies
\be\label{xxtm-1-ll-x1}
\nonumber z_3\in D\backslash
\mathbb{B}_{\delta_D}(z_1,3^{s-1}\delta_D(z_1,z_2))\;\;\mbox{and}\;\;
z_3\in \mathbb{B}_{\delta_D}(z_1,3^s\delta_D(z_1,z_2)).
\ee
Obviously, $s\geq 2$.
It is possible
that $\delta_D(z_1,z_{3})= 3^{s-1}\delta_D(z_1,z_2)$.
We let $x_1=z_1$, and let $x_2$, $\dots$, $x_s$ be points
such that for each $i\in\{1$, $\dots$, $s-1\}$, $x_{i+1}$ denotes the
last point in $\beta$ along the direction from $z_1$ to $z_3$ such
that (see Figure \ref{hconj-fig01})
\be\label{wang-laowu-1}
\delta_D(x_1,x_{i+1})= 3^{i-1}\delta_D(z_1,z_2)
.
\ee

\noindent
Apparently,
\beq \label{sh-1}
3^{s-1}\delta_D(z_1,z_2)\leq
\delta_D(z_1,z_3)<3^s\delta_D(z_1,z_2).
\eeq
Hence, we have chosen points $\{x_1, \ldots, x_s\}$ from $\beta$ including the
case $s=1$, and now, we still need to pick up another point, denoted by $y_{s+1}$,
in $\beta$, which is the first point in $\beta[z_3,x_s]$ along the
direction from $z_3$ to $x_s$ such that (see Figure \ref{hconj-fig01})
\be\label{wang-laowu-2}
\delta_D(y_{s+1},z_3)=\frac{1}{8}\delta_D(x_s,z_3).
\ee
Then we see that
\beq\label{srwh-1}
\delta_D(x_s,y_{s+1})\leq
\delta_D(x_s,z_3)+\delta_D(y_{s+1},z_3) =9\delta_D(y_{s+1},z_3),
\eeq
and for each $w\in \beta[y_{s+1},z_3]$,
\beq\label{plm-1}
\delta_D(y_{s+1},w) &\leq &
\delta_D(y_{s+1},z_3)+\delta_D(z_3,w) \leq \delta_D(y_{s+1},z_3)+
\frac{1}{8}\delta_D(x_s,z_3)\\ \nonumber
&\leq & 2\delta_D(y_{s+1},z_3).
\eeq
Moreover,
$$\delta_D(x_1,z_3)\geq \frac{3}{4}\delta_D(x_s,z_3).
$$

This inequality is obvious if $s=1$, and if $s\geq 2$, \eqref{wang-laowu-1} and (\ref{sh-1}) imply
$$\delta_D(x_1,z_3)\geq \delta_D(x_s,z_3)-\delta_D(x_1,x_s)\geq
\delta_D(x_s,z_3)-\frac{1}{3}\delta_D(x_1,z_3),
$$
from which the inequality easily follows.

Next, we have the following useful inequalities related to $y_{s+1}$. First,
we deduce from (\ref{wang-laowu-2}), (\ref{plm-1}) and the choice of $y_{s+1}$ that for
$w\in \beta[y_{s+1},z_3]$,
\beq\label{srwh-1-1}
\delta_D(x_1,w)&\geq&\delta_D(x_1,z_3)-\delta_D(w,z_3) \geq
\frac{3}{4}\delta_D(x_s,z_3)-\delta_D(w,z_3)
\\ \nonumber&\geq& 5\delta_D(y_{s+1},z_3)\\ \nonumber&\geq&\frac{5}{2}\delta_D(y_{s+1},w),
\eeq
and then, if $s\geq 2$, we infer from \eqref{sh-1} and (\ref{wang-laowu-2}) that
\beq\label{srwh-0}
\delta_D(x_s,y_{s+1})&\leq&
\delta_D(x_1,x_s)+\delta_D(x_1,z_3)+\delta_D(y_{s+1},z_3)
\\ \nonumber&=& \delta_D(x_1,x_s)+\delta_D(x_1,z_3)+\frac{1}{8}\delta_D(x_s,z_3)
\\ \nonumber&\leq& \frac{9}{8}\big(\delta_D(x_1,x_s)+\delta_D(x_1,z_3)\big)
\\ \nonumber&<&\frac{45}{4}\delta_D(x_1,x_s).
\eeq

We now take $v_i=x_i$ for each $i\in \{1,\ldots,s\}$, where $s\geq 1$,
$v_{s+1}=y_{s+1}$ and $v_{s+2}=z_3$. In this way, we get the desired partition of
$\beta$ (see Figure \ref{hconj-fig01}).

The following two claims are inequalities on the image of the partition of $\beta$ under $f$.
\bcl\label{xtwmm-1} If $s=1$, then $\delta_{D'}(v'_1,v'_3)\leq \mu_6\delta_{D'}(v'_1,z'_2)$. \ecl

Suppose on the contrary that
$$\delta_{D'}(v'_1,v'_3)>\mu_6 \delta_{D'}(v'_1,z'_2).
$$
We shall get a contradiction by using the conformal modulus of families of curves (in the following, we briefly say ``by using the conformal modulus").

\begin{figure}[htbp]
\begin{center}
\input{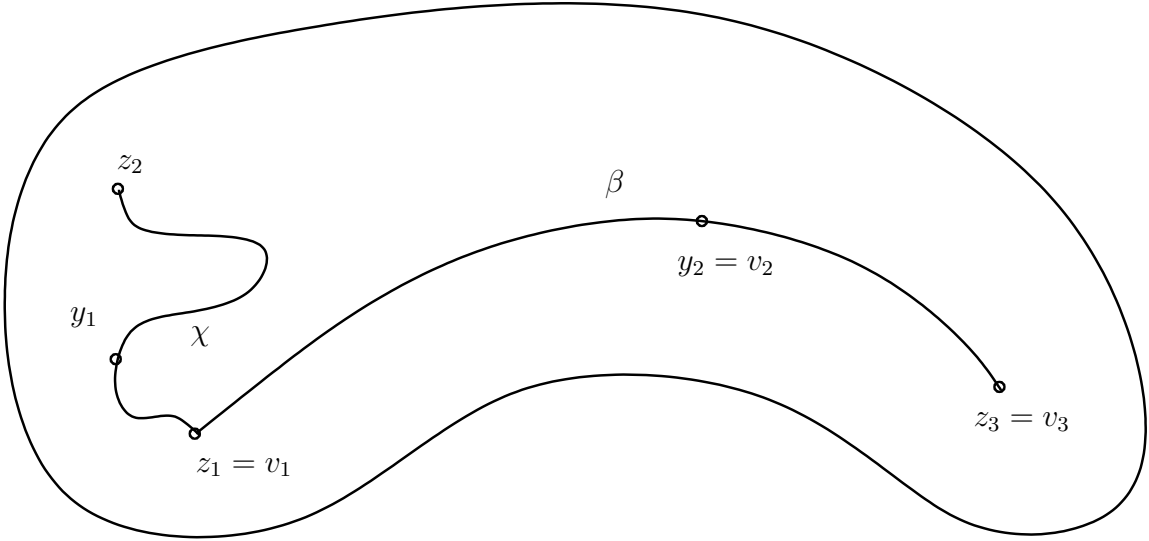}

\caption{The arc $\chi$ and the point $y_1$.}

\label{fig-0002}
\end{center}
\end{figure}

Let $\chi$ be an arc joining $v_1$ and $z_2$ in $A$ (see Figure \ref{fig-0002}), and let $y_1$
be the first point in $\chi$ along the direction from $v_1$ to $z_2$ satisfying
(see Figure \ref{fig-0002})
\be\label{samy-wang-1}
\delta_D(v_1,y_1)=\frac{1}{2}\delta_D(v_1,z_2).
\ee

First, we apply the assumption ``$D$ being $\varphi$-broad" to get a lower bound for the conformal modulus ${\rm Mod}(\chi[v_1,y_1],
\beta[v_2,v_3]; D)$. For this, we need to show that the curves $\chi[v_1,y_1]$ and $\beta[v_2,v_3]$ are disjoint.

For each $u\in\chi[v_1,y_1]$, by the assumption on the weak quasisymmetry of $f$ and the choice of $y_1$, we have
\be\label{samy-wang-2} \delta_{D'}(v'_1,u')\leq H\delta_{D'}(v'_1,z'_2),\ee
and for each $v\in\beta[v_2,v_3]$, it follows from (\ref{wang-laowu-2}) together with the choice of $y_1$
and $v_2$ that
\beq\label{whantti-1}
\delta_D(u,v)&\geq&
\delta_D(v_1,v_3)-\delta_D(v_3,v)-\delta_D(v_1,u)\\ \nonumber&\geq& \frac{7}{8}\delta_D(v_1,v_3)
-\frac{1}{2}\delta_D(v_1,z_2)\\ \nonumber&\geq&\frac{3}{8}\delta_D(v_1,v_3)= 3\delta_D(v_2,v_3)
\\ \nonumber&\geq& 3\delta_D(v,v_3).
\eeq So the assumption on the weak quasisymmetry of $f$ implies
\beq\label{antti-wang-0}
\delta_{D'}(v'_1,v'_3)&\leq&
\delta_{D'}(v'_1,v')+\delta_{D'}(v',v'_3)\leq
(H+1)\delta_{D'}(v'_1,v').
\eeq
Now, (\ref{samy-wang-2}) and (\ref{antti-wang-0}) show that
\beq\label{wang-manzi-3}
\delta_{D'}(u',v')&\geq&
\delta_{D'}(v'_1,v')-\delta_{D'}(v'_1,u')\\
\nonumber&\geq& \frac{1}{H+1}\delta_{D'}(v'_{1},v'_3)-H\delta_{D'}(v'_1, z'_2)
\\
\nonumber&>&
 \Big (\frac{\mu_6}{H+1}-H \Big )\delta_{D'}(v'_{1},z'_2).
\eeq
Hence $\chi[v_1,y_1]\cap \beta[v_2,v_3]=\emptyset$. Further, the combination of
(\ref{wang-laowu-2}) and (\ref{samy-wang-1}) guarantees that
$$\frac{\delta_D(v_1,v_3)}{\min\{\diam(\chi[v_1,y_1]),\diam(\beta[v_2,v_3])\}}\leq 8,
$$
we see from the assumption ``$D$ being $\varphi$-broad"
that
\beq
\nonumber {\rm Mod}(\chi[v_1,y_1],
\beta[v_2,v_3]; D)
&\geq& \varphi\Big(\frac{\delta_D(\chi[v_1,y_1], \beta[v_2,v_3])
}{\min\{\diam(\chi[v_1,y_1]),\diam(\beta[v_2,v_3])\}}\Big)\\ \nonumber
&\geq& \varphi\Big(\frac{\delta_D(v_1,v_3)}{\min\{\diam(\chi[v_1,y_1]),\diam(\beta[v_2,v_3])\}}\Big)
\\ \nonumber&\geq&\varphi(8),
\eeq where
$\varrho_n$ is from Theorem \Ref{lemm1--0''}, which is the desired bound.

Since \eqref{samy-wang-2} implies
$$ \diam(\chi'[v'_1,y'_1])\leq 2H\delta_{D'}(v'_1,z'_2)
$$
and the quasiconformal invariance property of the moduli of the families of curves shows
$${\rm Mod}(\chi[v_1,y_1],
\beta[v_2,v_3]; D) \leq K{\rm Mod}(\chi'[v'_1,y'_1],
\beta'[v'_2,v'_3]; D'),
$$
we infer from Theorem \Ref{lemm1--0''} and (\ref{wang-manzi-3}) that
\begin{eqnarray*}
\varphi(8)&\leq& {\rm Mod}(\chi[v_1,y_1],
\beta[v_2,v_3]; D) \leq K{\rm Mod}(\chi'[v'_1,y'_1],
\beta'[v'_2,v'_3]; D')\\
\nonumber&\leq& K \varrho_n\left(
\frac{(\frac{\mu_6}{H+1}-H)\delta_{D'}(v'_1,z'_2)}{\diam(\chi'[v'_1,y'_1])}\right)
\leq K \varrho_n\left(
\frac{\frac{\mu_6}{H+1}-H}{2H}\right)\\ \nonumber&<&
\varphi(8).
\end{eqnarray*}
This obvious contradiction shows that the claim is true.

\bcl\label{xtwmm-1-1}
If $s\geq 2$, then $\delta_{D'}(v'_1,v'_{i+1})\leq \mu_6\delta_{D'}(v'_1,v'_{i})$ for each $i\in \{2, \ldots, s+1\}$. \ecl

Suppose on the contrary
that there exists an $i\in \{2, \ldots, s+1\}$ such that
$$\delta_{D'}(v'_1,v'_{i+1})> \mu_6\delta_{D'}(v'_1,v'_{i}).
$$

If $i=s+1$, then we get
\begin{eqnarray*}
\delta_{D'}(v'_{s+1},v'_{s+2})&\geq& \delta_{D'}(v'_1,v'_{s+2})-
\delta_{D'}(v'_1,v'_{s+1})\\
&\geq& (\mu_6-1)\delta_{D'}(v'_{1},v'_{s+1}),
\end{eqnarray*}
whence, again, the assumption on the weak quasisymmetry of $f$ implies
$$\delta_D(v_{s+1},v_{s+2})\geq \delta_D(v_1,v_{s+1}),
$$
which contradicts (\ref{srwh-1-1}). Hence we see that $i\in\{2, \ldots, s\}$.
In the following, by using the conformal modulus, we shall show that it is impossible either.
For this purpose, some preparation is needed.

\begin{figure}[htbp]
\begin{center}
\input{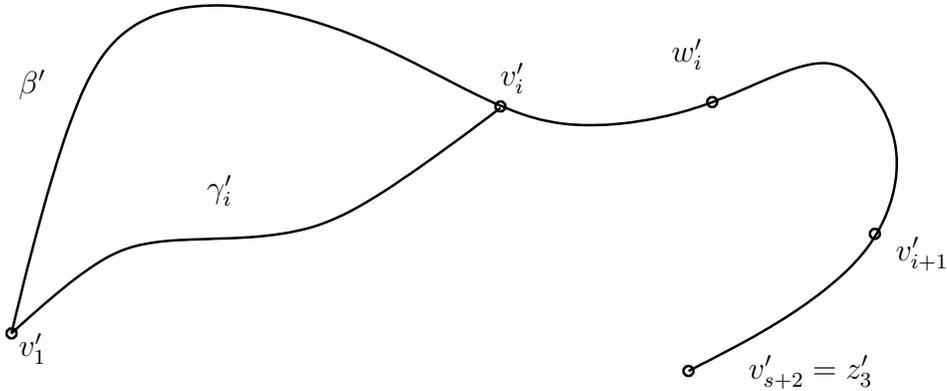}
\caption{The point $w'_i$, and the arcs $\beta'$ and $\gamma'_i$ in
$D'$}

\label{hconj-fig02}
\end{center}
\end{figure}

Let $\gamma'_i$ be an arc joining $v'_1$ and $v'_i$ in $D'$
such that (see Figure \ref{hconj-fig02})
\be\label{srwh-0-1}
\diam(\gamma'_i)<\frac{5}{4}\delta_{D'}(v'_1,v'_i), \ee and we let
$w'_i$ be the last point in $\beta'[v'_i,v'_{s+2}]$
along the direction from $v'_i$ to $v'_{s+2}$ such that (see Figure \ref{hconj-fig02})
\be\label{wang-laowu-3}
\delta_{D'}(v'_1,w'_i)=\frac{\mu_6}{3H}\delta_{D'}(v'_1,v'_i).
\ee
Obviously,
\be\label{srwh-0-3}
\delta_{D'}(v'_1,v'_{i+1})> 3H\delta_{D'}(v'_1,w'_i).
\ee

Let us leave the proof of the claim for a moment and determine the position
of $w'_i$ in $\beta'[v'_i,v'_{s+2}]$. Our result is as follows.
\bscl
\label{srwh-0-2} $w'_i\in \beta'[v'_i,v'_{i+1}].$
\escl

To prove this subclaim, we divide the arguments into the case where
$i\leq s-1$ and the case where $i=s$. We first consider the case
$i\leq s-1$. Then, once again, by the assumption on the weak quasisymmetry of $f$ together with (\ref{srwh-0-3}), we have
$$\delta_D(v_1,v_{i+1})>\delta_D(v_1,w_i),
$$
and so  the choice of $v_{i+1}$ shows that $w'_i\in \beta'[v'_i,v'_{i+1}]$.

On the other hand, for the case $i=s$, (\ref{wang-laowu-3}) leads to
\begin{eqnarray*}
\delta_{D'}(v'_{s+1},w'_s)&\geq&
\delta_{D'}(v'_s,v'_{s+1})-
\delta_{D'}(v'_1,v'_s)-\delta_{D'}(v'_1,w'_s)\\
&\geq& \delta_{D'}(v'_{1},v'_{s+1})-
2\delta_{D'}(v'_1,v'_s)-\delta_{D'}(v'_1, w'_s)
\\
\nonumber&\geq& \Big(\mu_6-2-\frac{\mu_6}{3H}\Big)
\delta_{D'}(v'_1,v'_s)
\\
\nonumber&>&
H\delta_{D'}(v'_1,w'_s),
\end{eqnarray*}
whence we infer from the assumption on the weak quasisymmetry of $f$ that
$$\delta_D(v_{s+1},w_s)> \delta_D(v_1,w_s),
$$
which together with (\ref{srwh-1-1}) shows that $w'_s\in\beta'[v'_s,v'_{s+1}]$.
Hence Subclaim \ref{srwh-0-2} is also true in this case. The proof of the subclaim is complete.
\medskip

We shall now present a proof of Claim \ref{xtwmm-1-1}. We shall reach a contradiction by obtaining
a lower bound and an upper bound for the conformal modulus ${\rm Mod}(\gamma_i,
\beta[w_i,v_{s+2}]; D)$. For a lower bound of this quantity, we need an inequality:
If $s\geq 2$, then for each $i\in\{2,\ldots,s\}$,
\beq\label{srwh-0-0}
\delta_D(v_i,v_{i+1})\leq
\frac{45}{4}\min\{\delta_D(v_1,v_i),\delta_D(v_{i+1},v_{s+2})\}.
\eeq

This inequality easily follows from (\ref{srwh-1}) and (\ref{srwh-0}) together with the following inequalities:
For $s\geq 3$ and $i\in
\{2,\ldots,s-1\}$,
$$\delta_D(v_i,v_{i+1})\leq
\delta_D(v_1,v_i)+\delta_D(v_1,v_{i+1})= 4\delta_D(v_1,v_i),
$$
and further
$$\delta_D(v_i,v_{i+1})\leq 4\delta_D(v_1,v_i)= \frac{4}{3}\delta_D(v_1,v_{i+1})\leq
\frac{2}{3}\delta_D(v_{i+1},v_{s+2}),
$$
where in the third inequality, the following estimate is used:
\begin{eqnarray*}
\delta_D(v_{i+1},v_{s+2}) & \geq &
\delta_D(v_{1},v_{s+2})-\delta_D(v_1,v_{i+1})
\geq 3^{s-1}\delta_D(z_{1},z_2)-\delta_D(v_1,v_{i+1})\\
&\geq & 2\delta_D(v_1,v_{i+1}).
\end{eqnarray*}

Since for each $w'\in\gamma'_i$ and $z'\in\beta'[w'_i,v'_{s+2}]$, we have
from \eqref{srwh-0-1} and the choice of $w_i'$ that
\beq\label{srwh-4}
\delta_{D'}(w',z')\geq
\delta_{D'}(v'_1,z')-\delta_{D'}(v'_1,w')\geq
\Big(\frac{\mu_6}{3H}-\frac{5}{4}\Big)\delta_{D'}(v'_1,v'_i),
\eeq
we see that $\gamma_i\cap \beta[w_i,v_{s+2}]=\emptyset$, whence it follows from
the assumption ``$D$ being $\varphi$-broad"  that
$$\ {\rm Mod}(\gamma_i,
\beta[w_i,v_{s+2}]; D)\geq\varphi\left(\frac{\delta_D(\gamma_i,
\beta[w_i,v_{s+2}])}{\min\{\diam(\gamma_i),\diam(\beta[w_i,v_{s+2}])\}}\right).
$$
Moreover, (\ref{srwh-0-0}) and Subclaim \ref{srwh-0-2} show that
\begin{eqnarray*}
\delta_D(\gamma_i, \beta[w_i,v_{s+2}])
\leq \delta_D(v_i,v_{i+1}) \leq
\frac{45}{4}\min\{\diam(\gamma_i),\diam(\beta[w_i,v_{s+2}])\},
\end{eqnarray*}
and so we obtain a lower bound as follows:
\beq\label{srwh-3}
{\rm Mod}(\gamma_i, \beta[w_i,v_{s+2}]; D)\geq\varphi\big(\frac{45}{4}\big).
\eeq

\begin{figure}[htbp]
\begin{center}
\input{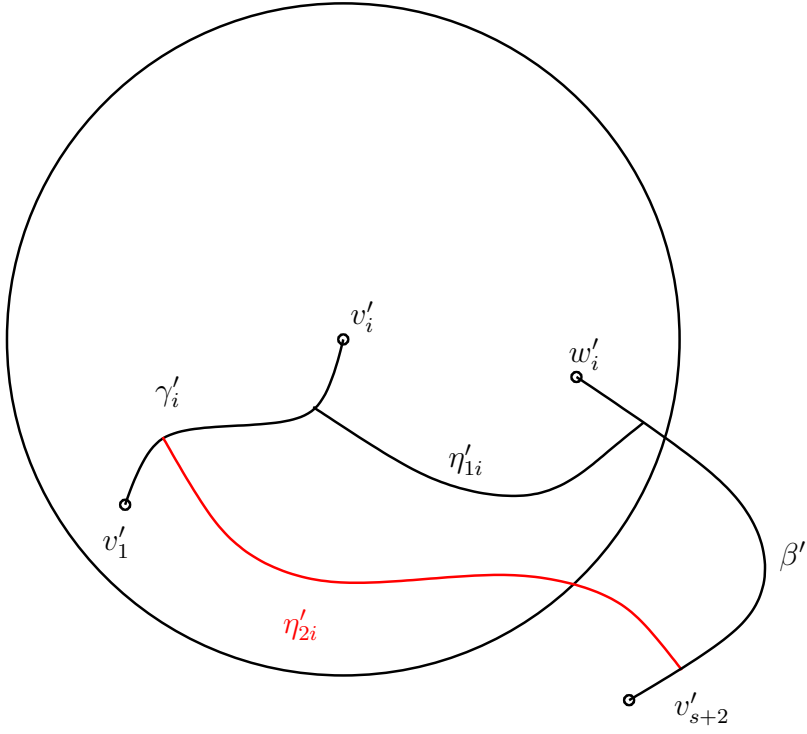}
\caption{The arcs $\eta'_{1i}\in \Upsilon'_{1i}$ and $\eta'_{2i}\in
\Upsilon'_{2i}$ in $D'$}

\label{hconj-fig03}
\end{center}
\end{figure}

For an upper bound for ${\rm Mod}(\gamma_i,
\beta[w_i,v_{s+2}]; D)$, we still need to construct two families of curves in $D'$ (see Figure \ref{hconj-fig03}):
$$\Upsilon'_{1i}=\{\eta'_{1i}: \eta'_{1i}\subset \mathbb{B}\big(v'_i, (\frac{\mu_6}{4H})^{\frac{1}{2}}
\delta_{D'}(v'_1,v'_i)\big)\}
$$
and
$$\Upsilon'_{2i}=\{\eta'_{2i}: \eta'_{2i}\cap  \mathbb{S}\big(v'_i, (\frac{\mu_6}{4H})^{\frac{1}{2}}
\delta_{D'}(v'_1,v'_i)\big)\not=\emptyset\},
$$
where $\eta'_{ji}$ $(j=1,2)$ denote the curves joining $\gamma'_i$
and $\beta'[w'_i,v'_{s+2}]$ in $D'$.

Now, the quasiconformal invariance property of moduli of the families of curves implies
$${\rm Mod}(\gamma_i,
\beta[w_i,v_{s+2}]; D)  \leq
K{\rm Mod}(\gamma'_i, \beta'[w'_i,v'_{s+2}]; D'),
$$
whence it follows from Theorems \Ref{lemm1--0''} and \Ref{lemm1--xx} together with
(\ref{srwh-4}) that
\beq\label{wlw-1}
\nonumber {\rm Mod}(\gamma_i, \beta[w_i,v_{s+2}]; D)
&\leq& K\big({\rm Mod}(\Upsilon'_{1i})+{\rm
Mod}(\Upsilon'_{2i})\big)\\
\nonumber&\leq& K\omega_{n-1}\left(
\frac{(\frac{\mu_6}{4H})^{\frac{n}{2}}}{(\frac{\mu_6}{3H}-\frac{5}{4})^n}
+\Big(\frac{1}{2}\log
\frac{4\mu_6}{25H}\Big)^{1-n}\right)\\
\nonumber&<&
\varphi\big(\frac{45}{4}\big),
\eeq
where in the second inequality, the fact used is $\eta'_{2i}\cap \mathbb{S}\big(v'_i,
\frac{5}{4}\delta_{D'}(v'_1,v'_i)\big)\not=\emptyset$ for each $\eta'_{2i}\in
\Upsilon'_{2i}$, which easily follows from
\eqref{srwh-0-1}. This obviously contradicts (\ref{srwh-3}). Hence Claim
\ref{xtwmm-1-1} is true.\medskip

Let us now finish the proof of the lemma. Let $\mu_7(c)=H\mu_6^{1+\log_3c}$. Then, since $\mu_6< \mu_7(c)$,
Claim \ref{xtwmm-1} guarantees that Lemma \ref{Lemc1--0-00} is true when $s=1$. For the case $s\geq 2$,
we see from the assumption on the weak quasisymmetry of $f$ and Claim \ref{xtwmm-1-1} that
$$\delta_{D'}(z'_1,z'_2)\geq \frac{1}{H}\delta_{D'}(v'_1,v'_2)\geq
\frac{1}{H\mu_6}\delta_{D'}(v'_1,v'_{3})\geq\ldots\geq
\frac{1}{H\mu_6^s}\delta_{D'}(v'_1,v'_{s+2}),
$$
and so
$$ \delta_{D'}(z'_1,z'_3)\leq H\mu_6^s \delta_{D'}(z'_1,z'_2)  \leq H\mu_6^{1+\log_3c}\delta_{D'}(z'_1,z'_2),
$$
since \eqref{sh-1} and the assumption ``$\delta_D(z_1,z_3)\leq c\delta_D(z_1,z_2)$" imply that $s-1\leq \log_3c$,
which shows that the lemma is also true in this case. Hence the proof of Lemma \ref{Lemc1--0-00} is complete.
\epf

\begin{lem}\label{mzz-1}
Under the assumptions of Lemma \ref{Lemc1--0-00}, we see that for all $a$, $x$, $y\in A$, if
$\delta_{D}(a,x)\leq \delta_{D}(a,y)$, then
$$\frac{\delta_{D'}(a',x')}{\delta_{D'}(a',y')}\leq \psi\Big(\frac{\delta_{D}(a,x)}{\delta_{D}(a,y)}\Big),
$$
where $\psi:\,(0,1]\to (0, \infty)$ is an
increasing homeomorphism which depends only on the data
$$\sigma=(n,K,H,\varphi, \varrho_n).
$$
\end{lem}

\bpf For a proof, we let $a$, $x$ and $y\in A$ with $\delta_{D}(a,x)\leq \delta_{D}(a,y)$.
For convenience, we write
$$r=\frac{\delta_D(a,y)}{\delta_D(a,x)}\;\;\mbox{and}\;\;s=\frac{\delta_{D'}(a',y')}{\delta_{D'}(a',x')}.
$$
Obviously, $r\geq 1$. Let
$$\ds \mu_8=\max\left\{6^3K, H^2\exp\left ({\Big(\frac{8K\omega_{n-1}}{\varphi(6)}\Big)^{\frac{1}{n-1}}}\right)\right \}.
$$
With the aid of $\mu_8$, we divide the discussions into two cases: $r\leq 4\mu_8^2$ and  $r> 4\mu_8^2$.
For each case, we shall construct
an increasing homeomorphism. The desired homeomorphism follows from this.

Suppose first that $r\leq 4\mu_8^2$. Then by the assumption ``$f|_A$ being
weakly $H$-QS", we have
$$\frac{\delta_{D'}(a',x')}{\delta_{D'}(a',y')}\leq H\leq4H\mu_8^2\frac{\delta_D(a,x)}{\delta_D(a,y)}.
$$
In this case we let
\beq\label{tue-2}
\psi_1(p)=4H\mu_8^2p
\eeq
for $p$ in $[\frac{1}{4\mu_8^2},1]$.

\begin{figure}[htbp]
\begin{center}
\input{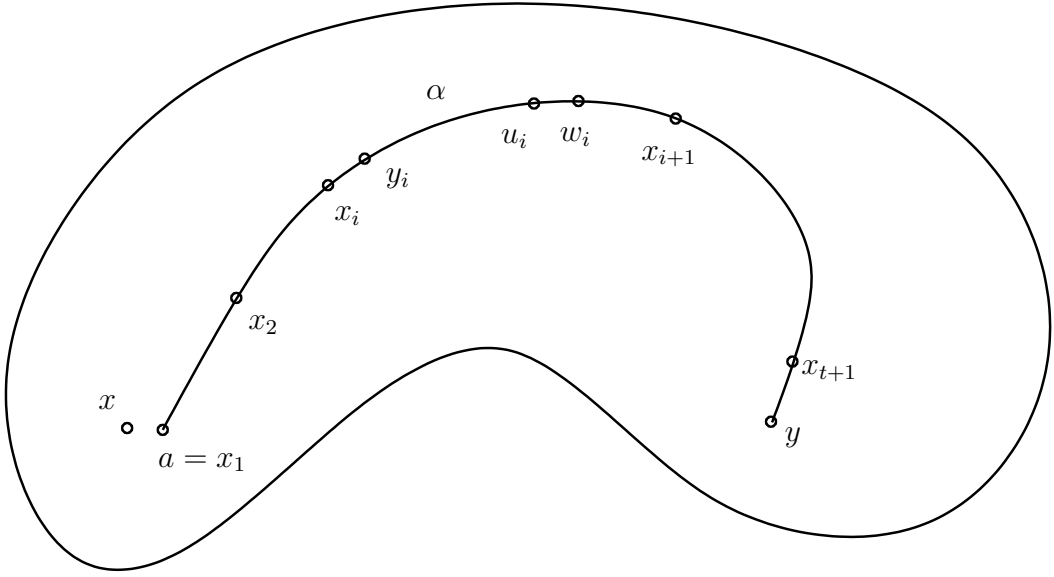}
\caption{The arc $\alpha$ and the related points in $A$ 
}

\label{hconj-fig04}
\end{center}
\end{figure}

Suppose next that $r> 4\mu_8^2$. In this case, we shall exploit the conformal modulus together with Theorems \Ref{lemm1--0''} and \Ref{lemm1--xx} to get an increasing homeomorphism. First, we do some preparation.

Let $\alpha$ be an arc
joining $a$ and $y$ in $A$ (see Figure \ref{hconj-fig04}). We give the following partition to $\alpha$. Let $x_1=a$, and let $x_2$, $\dots$,
$x_{t+1}$ be points such that for each $i\in\{1$, $\dots$, $t\}$,
$x_{i+1}$ denotes the last point in $\alpha$ along the direction
from $a$ to $y$ such that (see Figure \ref{hconj-fig04})
\be\label{xxtm-1}
\delta_D(x_i,x_{i+1})=
6^i\delta_D(a,x)\;\;\mbox{and}\;\;\delta_D(x_{t+1},y)<6\delta_D(x_t,x_{t+1}).
\ee
It is possible that $x_{t+1}=y$. Clearly, for this partition, we have
$$\sum_{i=2}^{t+1}\delta_D(x_{i-1}, x_i)=
\frac{6}{5}\Big(\delta_D(x_t,x_{t+1})-\delta_D(a,x)\Big),
$$
whence
$$\delta_D(a,y)\leq\sum_{i=2}^{t+1}\delta_D(x_{i-1}, x_i)+
\delta_D(x_{t+1},y)< \frac{36}{5}\delta_D(x_t,x_{t+1}),
$$
and so
$$\sum_{i=2}^{t+1}\delta_D(x_{i-1}, x_i)\geq
\Big(\frac{1}{6}-\frac{3}{10\mu_8^2}\Big)\delta_D(a,y),
$$
since
$\frac{\delta_D(a,y)}{\delta_D(a,x)}>4\mu_8^2$. Thus we infer from
$$\sum_{i=2}^{t+1}\delta_D(x_{i-1}, x_i)=\frac{6}{5}(6^t-1)\delta_D(a,x)
$$
that
\be\label{mxt-2-5}
t> \log_{6} (4r)-2.
\ee

For each $i\in\{2,\ldots,t\}$, we still need to pick up three points from $\alpha[x_i, x_{i+1}]$ as follows. Let $y_i$ (resp. $u_i$) denote the
first point in $\alpha$ along the direction from $x_i$ to $x_{i+1}$
(resp. from $y_i$ to $x_{i+1}$) such that (see Figure \ref{hconj-fig04})
\be\label{mm-1}
\delta_D(x_i,y_i)=\frac{1}{\mu_9^2}\delta_D(x_{i-1},x_i)\;\;
\big(\mbox{resp.}\;
\delta_D(x_i,u_i)=\frac{6}{\mu_9^2}\delta_D(x_{i-1},x_i)\big),
\ee
and let $w_i$  be the first point in $\alpha$ along the direction from $u_i$ to $x_{i+1}$
such that (see Figure \ref{hconj-fig04})
\be\label{mm-2}
\delta_D(u_i,w_i)=\frac{1}{\mu_9^2}\delta_D(x_{i-1},x_i).
\ee
We see that
$$\alpha[x_i,y_i]\cap \alpha[u_i,w_i]=\emptyset.
$$
This can be seen from the following estimate: For $w\in \alpha[x_i,y_i]$ and $z\in\alpha[u_i,w_i]$, the choice
of $y_i$ and $w_i$ implies that
\beq\label{mm-4}
\delta_D(w,z)&\geq&
\delta_D(x_i,u_i)-\delta_D(u_i,z)-\delta_D(x_i,w)\\
\nonumber&\geq&
\frac{4}{\mu_9^2}\delta_D(x_{i-1},x_i). \eeq
Thus the assumption the ``$D$ being $\varphi$-broad" together with (\ref{mm-1}) and (\ref{mm-2}) leads to
\beq\label{mm-3}
{\rm Mod}(\alpha[x_i,y_i], \alpha[u_i,w_i]; D)&\geq&
\varphi\Big(\frac{\delta_D(\alpha[x_i,y_i],
\alpha[u_i,w_i])}{t_i}\Big)\\
\nonumber&\geq&\varphi\Big(\frac{\delta_D(x_i,u_i)}{t_i}\Big)\\
\nonumber&\geq&\varphi(6),
\eeq where
$t_i=\min\{\diam(\alpha[x_i,y_i]),\diam(\alpha[u_i,w_i])\}$.

\begin{figure}[htbp]
\begin{center}
\input{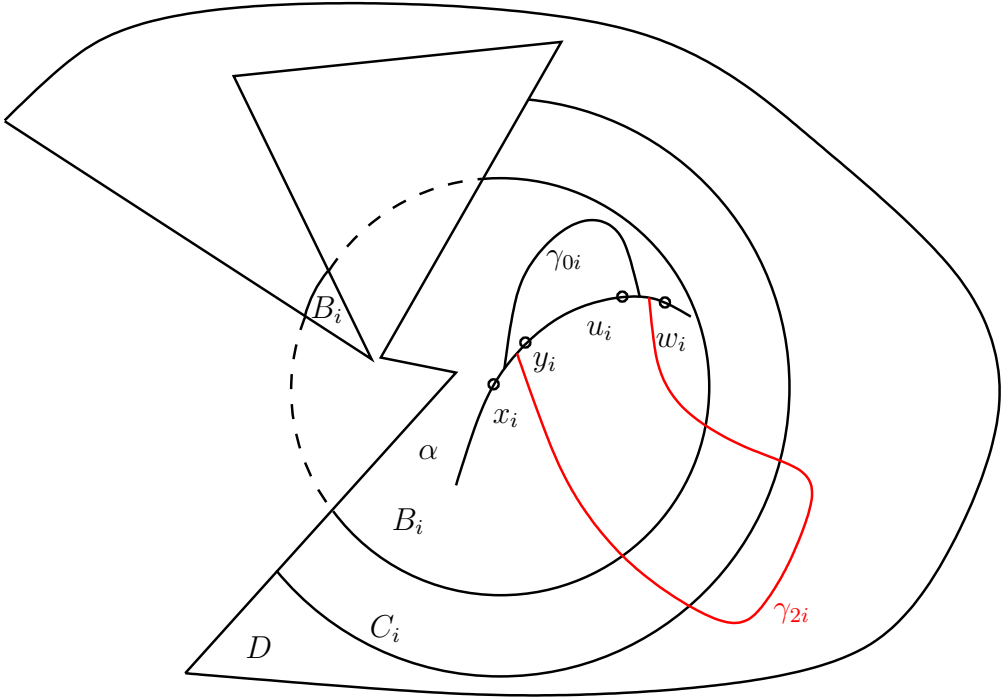}
\caption{The arcs $\gamma_{0i}\in\Gamma_{0i}$ and
$\gamma_{2i}\in\Gamma_{2i}$ in $D$}

\label{hconj-fig05}
\end{center}
\end{figure}

Let
$\Gamma_i$ denote the family of all curves connecting $\alpha[x_i, y_i]$ and $\alpha[u_i, w_i]$ in $D$. Then
$${\rm Mod}(\Gamma_i)={\rm Mod}(\alpha[x_i,y_i], \alpha[u_i,w_i]; D).
$$
To decompose $\Gamma_i$, we construct a finite sequence of balls in $D$ as follows.
For each $i\in\{2,\dots,t\}$, we let
$$C_i=\mathbb{B}_{\delta_D}\Big(x_i,\frac{1}{\mu_9}\delta_D(x_{i-1},x_i)\Big),
$$
where $\mu_9=\exp({4\mu_8})$. Then (see Figure \ref{hconj-fig05})
$$\Gamma_i=\Gamma_{0i}\cup \Gamma_{1i}\cup \Gamma_{2i}
$$
where
\beqq
\Gamma_{0i}&=&\{\gamma\in \Gamma_i: \gamma\subset  B_i\cap C_i\},\\
\Gamma_{1i}&=&\{\gamma\in \Gamma_i: \gamma\subset B_i-B_i\cap C_i\},\\
\Gamma_{2i}&=&\left \{\gamma\in \Gamma_i: \gamma\cap \mathbb{S}\Big(x_i,
\frac{\mu_8}{\mu_9^2}\delta_D(x_{i-1},x_i)\Big)\not=\emptyset\right\},\\
B_i&=&\mathbb{B}\Big(x_i, \frac{\mu_8}{\mu_9^2}\delta_D(x_{i-1},x_i)\Big).
\eeqq
Thus, we see that
\bcl\label{in-fin-1} $\Gamma_{1i}=\emptyset \;\;\mbox{(see Figure \ref{hconj-fig05})}.$
\ecl

For a proof of this claim, it suffices to show that for each $\gamma$ in $\Gamma_i$, $\gamma\subset B_i$ implies
$\gamma\subset C_i$. Obviously, for each $\gamma$ in $B_i$, we have
$$\diam(\gamma)\leq
\frac{2\mu_8}{\mu_9^2}\delta_D(x_{i-1},x_i),
$$
and so by (\ref{mm-1}),
$$\delta_D(x_i,x)\leq \delta_D(x_i,y_i)+\diam(\gamma)\leq
\Big(\frac{1+2\mu_8}{\mu_9^2}\Big)\delta_D(x_{i-1},x_i)<\frac{1}{\mu_9}\delta_D(x_{i-1},x_i)
$$
for each $x\in \gamma$, which implies $\gamma\subset C_i$, and so Claim \ref{in-fin-1} is proved.

\medskip
It follows from Claim \ref{in-fin-1} that
$$\Gamma_i=\Gamma_{0i}\cup \Gamma_{2i},
$$
and then the choice of $y_i$, (\ref{mm-3}) and Theorem \Ref{lemm1--xx} imply
\begin{eqnarray*}
\varphi(6)&\leq&{\rm Mod}(\Gamma_i)
\leq {\rm Mod}(\Gamma_{0i})+{\rm Mod}(\Gamma_{2i})  \leq {\rm Mod}(\Gamma_{0i})+\omega_{n-1}(\log
\mu_8)^{1-n}\\
\nonumber&\leq& {\rm Mod}(\Gamma_{0i})+\frac{1}{4}\varphi(6),
\end{eqnarray*}
where in the third inequality, we have used the fact that for each $\gamma\in \Gamma_{2i}$,
$$\gamma\cap \mathbb{S}\big(x_i, \frac{1}{\mu_9^2}\delta_D(x_{i-1},x_i)\big)\not=\emptyset .
$$
Now, the quasiconformal invariance property of moduli of the families of curves implies
\be\label{mm-5}
{\rm Mod}(\Gamma'_{0i})\geq \frac{1}{K}{\rm Mod}(\Gamma_{0i})\geq \frac{3}{4K}\varphi(6),
\ee
which will be needed later on.

The following discussion still needs a lower bound on the inner distance between $\alpha'[x'_i,y'_i]$ and $\alpha'[u'_i,w'_i]$. We
first prove some elementary inequalities.
It follows from (\ref{llvv-1}) and (\ref{mm-1}) that
$$\delta_D(x_i,y_i)=\frac{1}{\mu_9^2}\delta_D(x_{i-1},x_i)\geq \frac{5}{6\mu_9^2}\delta_D(x_1,x_i)
\geq \frac{5}{\mu_9^2}\delta_D(a,x),
$$
whence Lemma \ref{Lemc1--0-00} implies
\be\label{xxllaa-1}
\delta_{D'}(x'_1,x'_i)\leq
\mu_7\big(\frac{6}{5}\mu_9^2\big)\delta_{D'}(x'_i,y'_i),
\ee
and the assumption on the weak quasisymmetry of $f$ shows
\be\label{xxllaa-1'}
\delta_{D'}(a',x')\leq H\delta_{D'}(x'_1,x'_i).
\ee

Now, we are ready to establish the lower bound.
For all $w\in \alpha[x_i,y_i]$ and $z\in\alpha[u_i,w_i]$, since
$$\delta_D(x_i,w)\leq \frac{1}{\mu_9^2}\delta_D(x_{i-1},x_i)
$$
and
$$\delta_D(w,y_i)\leq \delta_D(x_i,w)+\delta_D(x_i,y_i)
\leq\frac{2}{\mu_9^2}\delta_D(x_{i-1},x_i),
$$
we infer from (\ref{mm-4}) that
$$\delta_D(w,z)>\max\{\delta_D(x_i,w), \delta_D(w,y_i)\}.
$$
Then (\ref{xxllaa-1}) and \eqref{xxllaa-1'} show
\beq\label{mmvvxx-1}
\delta_{D'}(w',z')&\geq& \frac{1}{H}\max\{\delta_{D'}(x'_i,w'),\delta_{D'}(w',y'_i)\}\geq
 \frac{1}{2H}\delta_{D'}(x'_i,y'_i)\\ \nonumber&\geq&
\frac{1}{2H^2\mu_7\big(\frac{6}{5}\mu_9^2\big)}\delta_{D'}(a',x'),
\eeq
which is the required lower bound.

\begin{figure}[htbp]
\begin{center}
\input{figure06.pspdftex}
\caption{The arcs $\gamma'_{0i_1}\in\Gamma'_{0i_1}$ and
$\gamma'_{0i_2}\in\Gamma'_{0i_2}$ in $D'$
\label{hconj-fig06}}
\end{center}
\end{figure}

In order to apply Theorems \Ref{lemm1--0''} and \Ref{lemm1--xx} in the proof,
we decompose $\Gamma'_{0i}$ in the following way (see Figure \ref{hconj-fig06}):
$$\Gamma'_{0i}=\Gamma'_{0i_1}\cup \Gamma'_{0i_2} 
,
$$
where
$$\Gamma'_{0i_1}=\{\gamma'\in \Gamma'_{0i}: \;\gamma'\subset
\mathbb{B}(x'_i,\mu_8\delta_{D'}(a',y'))\}
$$
and
$$\Gamma'_{0i_2}=\{\gamma'\in \Gamma'_{0i}:\;  \gamma' \cap
\mathbb{S}(x'_i,\mu_8\delta_{D'}(a',y'))\not=\emptyset\}.
$$
Set
$$B'_{1i}=B'_i\cap C'_i\cap \mathbb{B}(x'_i,\mu_8\delta_{D'}(a',y')).
$$

At present, we shall obtain a relationship between the curve family $\Gamma'_{0i_2}$
and the sphere $\mathbb{S}(x'_i,H^2\delta_{D'}(a',y'))$. Since (\ref{xxtm-11-1}) implies
\beq\label{plm2}
\delta_D(a,y)&\geq& \delta_D(x_t,y)-\delta_D(x_1,x_t)\geq
\delta_D(x_t,x_{t+1})-\delta_D(x_1,x_t)\\ \nonumber&\geq&\frac{4}{5}\delta_D(x_t, x_{t+1}),
\eeq
and so for each $i\in\{2,\ldots,t\}$ and $w\in\alpha[x_i,y_i]$, the choice of $y_i$,
(\ref{xxtm-11-1}) and (\ref{llvv-1}) guarantee
\beq\nonumber
\delta_D(a,y)&\geq& \frac{4}{5}\delta_D(x_i,x_{i+1})\geq 4\delta_D(x_1,x_i)
\geq \frac{16}{5}\delta_D(x_{i-1},x_i)=\frac{16}{5}\mu_9^2\delta_D(x_i,y_i)\\ \nonumber
&\geq & \frac{16}{5}\mu_9^2\delta_D(x_i,w).\eeq Hence by the assumption on the weak quasisymmetry of $f$, we have
$$\delta_{D'}(a',y')\geq \frac{1}{H}\delta_{D'}(a',x'_i)\geq \frac{1}{H^2}\delta_{D'}(x'_i,w'),
$$
which implies that
$$\alpha'[x'_i,y'_i]\subset \overline{\mathbb{B}}(x'_i, H^2\delta_{D'}(a',y')),
$$
and thus, we see that for each $\gamma'\in \Gamma'_{0i_2}$ (see Figure \ref{hconj-fig06}),
$$\gamma'\cap \mathbb{S}(x'_i,H^2\delta_{D'}(a',y'))\not=\emptyset
.
$$

Now, we are ready to apply Theorems \Ref{lemm1--0''} and \Ref{lemm1--xx} to get an increasing homeomorphism.
It follows from (\ref{mm-5}) and Theorem \Ref{lemm1--xx} that
\begin{eqnarray*}
\frac{3}{4K}\varphi(6)&\leq& {\rm Mod}(\Gamma'_{0i})
\leq {\rm Mod}(\Gamma'_{0i_1})+{\rm Mod}(\Gamma'_{0i_2})\\
\nonumber&\leq& {\rm Mod}(\Gamma'_{0i_1})+\omega_{n-1}(\log
\frac{\mu_8}{H^2})^{1-n}\\
\nonumber&\leq& {\rm Mod}(\Gamma'_{0i_1})+\frac{1}{8K}\varphi(6),
\end{eqnarray*}
whence
$${\rm Mod}(\Gamma'_{0i_1})\geq \frac{5}{8K}\varphi(6),
$$
and so by (\ref{mmvvxx-1}) and Theorem \Ref{lemm1--0''},
$$\frac{5}{8K}\varphi(6)\leq {\rm Mod}(\Gamma'_{0i_1})\leq
\bigg(\frac{2H^2\mu_7\big(\frac{6}{5}\mu_9^2\big)}{\delta_{D'}(a',x')}\bigg)^n m(B'_{1i}),
$$
from which we get
\be\label{llxxm-1}
m(B'_{1i}) \geq
\frac{5\varphi(6)}{8K\bigg(2H^2\mu_7\big(\frac{6}{5}\mu_9^2\big)\bigg)^n}(\delta_{D'}(a',x'))^n.
\ee

Since we see from (\ref{xxtm-11-1}) and (\ref{plm2}) that for each $i\in\{2,\ldots,t\}$,
\begin{eqnarray*}
\delta_D(x_1,x_i)\leq\delta_D(x_1,x_{i-1})+\delta_D(x_{i-1},x_i)
\leq \frac{6}{5}\delta_D(x_{i-1},x_i) \leq \frac{1}{4}\delta_D(a,y),
\end{eqnarray*} so
again, the assumption on the weak quasisymmetry of $f$ implies that for each $w'\in \overline{B}'_{1i}$,
\begin{eqnarray*}
\delta_{D'}(x'_1,w')&\leq&\delta_{D'}(x'_1,x'_i)+\delta_{D'}(x'_i,w')
\leq (H+\mu_8)\delta_{D'}(a',y')< 2\mu_8\delta_{D'}(a',y'),
\end{eqnarray*}
which assures the inclusion
$$B'_{1i}\subset \mathbb{B}(x'_1,2\mu_8H\delta_{D'}(a',y')).
$$

The disjointness of the balls $\{C_i\}_{i=2}^t$ is needed now, which is indicated in the following claim.
\bcl
\label{xxll-1} For any $i\not= j\in\{2,\dots,t\}$,
$C_i\cap C_j=\emptyset.$
\ecl

It follows from (\ref{xxtm-1})
that for each $q<i\in\{1,\dots,t\}$,
\be\label{xxtm-11-1}
\delta_D(x_q, x_i)\leq\sum_{j=q}^{i-1}\delta_D(x_{j}, x_{j+1})=\frac{6}{5}(6^{i-q}-1)\delta_D(a,x)
< \frac{1}{5}\delta_D(x_i,x_{i+1}).
\ee
Thus we have that for $i\geq 3$ and $j\in\{1,\dots,i-2\}$,
\beq\label{llvv-1}
\frac{6}{5}\delta_D(x_{i-1}, x_i)&\geq&\delta_D(x_{i-1}, x_i)+\delta_D(x_j, x_{i-1})\\
\nonumber&\geq&\delta_D(x_j, x_i)\\
\nonumber&\geq& \delta_D(x_{i-1}, x_i)-\delta_D(x_j, x_{i-1})\\
\nonumber&\geq& \frac{4}{5}\delta_D(x_{i-1}, x_i),
\eeq
and so for all $i\not=j\in\{2,\ldots,
t\}$, it follows that
\begin{eqnarray*}
\delta_D(x_i,x_j)\geq \frac{4}{5}\max\{\delta_D(x_{i-1},
x_i),\delta_D(x_{j-1}, x_j)\}\geq
\frac{2}{5}\mu_9\max\{\delta_D(C_i),\delta_D(C_j)\},
\end{eqnarray*}
from which the claim follows.\medskip

We see from Claim \ref{xxll-1} and (\ref{llxxm-1}) that
\begin{eqnarray*}
\frac{5\varphi(6)}{8K\bigg(2H^2\mu_7\big(\frac{6}{5}\mu_9^2\big)\bigg)^n}(\delta_{D'}(a',x'))^n
\cdot t&\leq& \sum_{i=1}^{t}m(B'_{1i})
\\ \nonumber&\leq& m(\mathbb{B}(a',2\mu_8H\delta_{D'}(a',y')))\\
\nonumber&\leq&
\omega_{n-1}(2\mu_8H\delta_{D'}(a',y'))^n.
\end{eqnarray*} Hence
 (\ref{mxt-2-5}) leads to
$$\frac{\delta_{D'}(a',x')}{\delta_{D'}(a',y')}\leq
4H^2\mu_8\mu_7\big(\frac{6}{5}\mu_9^2\big)\Big(\frac{16K\omega_{n-1}}{5\varphi(6)\log_{6}
\frac{\delta_D(a,y)}{\delta_D(a,x)}}\Big)^{\frac{1}{n}}.
$$
In this case, we let
\beq\label{tue-1}
\psi_2(p)=4H^2\mu_8\mu_7\big(\frac{6}{5}\mu_9^2\big)
\Big(\frac{16K\omega_{n-1}}{5\varphi(6)\log_6\frac{1}{p}}\Big)^{\frac{1}{n}}
\eeq
for $p$ in $(0,\frac{1}{4\mu_8^2})$.

In conclusion, we see from \eqref{tue-2} and \eqref{tue-1} that the homeomorphism
$$\psi(p)=\begin{cases}
\displaystyle  \ds
\frac{\mu_{10}}{(\log_6\frac{1}{p})^{\frac{1}{n}}} &
\mbox{ if}\,\; \ds p\in(0,\frac{1}{4\mu_8^2}),\\
\displaystyle \ds \frac{4\mu_8^2\mu_{10}}{(\log_6
4\mu_8^2)^{\frac{1}{n}}} p  & \mbox{ if}\,\; \ds p\in[\frac{1}{4\mu_8^2},1]
\end{cases}
$$
is the desired one, where $\mu_{10}=\max\left \{4H\mu_8^2,
4H^2\mu_8\mu_7\big(\frac{6}{5}\mu_9^2\big)
\Big(\frac{16K\omega_{n-1}}{5\varphi(6)}\Big)^{\frac{1}{n}}\right \}$.
\epf

\subsection{The proof of Theorem \ref{thm1}}

To prove this theorem, by definition, it suffices to show that there exists an increasing homeomorphism
$\eta:\,(0,\infty)\to (0,\infty)$ such that the inequality
$$\frac{\delta_{D'}(x'_1,x'_2)}{\delta_{D'}(x'_1,x'_3)}\leq \eta\Big(\frac{\delta_D(x_1,x_2)}{\delta_D(x_1,x_3)}\Big)
$$
holds for all $x_1$, $x_2$ and $x_3\in A$ with $x_1\not= x_3$.
 We divide the construction into the cases where
$\delta_D(x_1,x_2)\leq\delta_D(x_1,x_3)$ and  where
$\delta_D(x_1,x_2)>\delta_D(x_1,x_3)$. In each case, we shall get a
homeomorphism or homeomorphisms. Then we construct the desired
homeomorphism $\eta$ from the obtained ones.

First, we suppose that $\delta_D(x_1,x_2)\leq\delta_D(x_1,x_3)$. Then Lemma \ref{mzz-1}
shows that
$$\frac{\delta_{D'}(x'_1,x'_2)}{\delta_{D'}(x'_1,x'_3)}\leq
\psi\Big(\frac{\delta_{D}(x_1,x_2)}{\delta_{D}(x_1,x_3)}\Big),
$$
where $\psi$ is the increasing homeomorphism constructed in Lemma \ref{mzz-1}. In
this case, we let
\beq\label{sun-1}
\eta_1(t)=\psi(t)
\eeq
for $t$ in $(0,1]$.

Next, we consider the case $\delta_D(x_1,x_2)>\delta_D(x_1,x_3)$.
Again, we divide the discussion into two cases which are as
follows.

\bca\label{xml-1} $\delta_{D'}(x'_1,x'_2)\leq
H_1\delta_{D'}(x'_1,x'_3)$, where $H_1=36H^2$.\eca

Apparently,
$$\frac{\delta_{D'}(x'_1,x'_2)}{\delta_{D'}(x'_1,x'_3)}\leq H_1\leq
H_1\frac{\delta_{D}(x_1,x_2)}{\delta_{D}(x_1,x_3)}.
$$
In this case, we define
\beq\label{fun-2}
\eta_2(t)=H_1\psi(1)t
\eeq
for $t> 1$.

\begin{figure}[htbp]
\begin{center}
\input{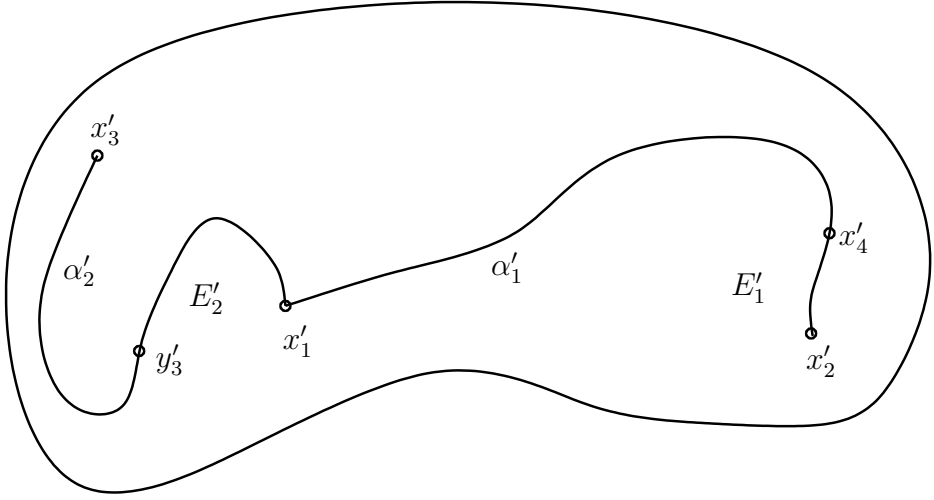}
\caption{The arcs $\alpha'_1$, $\alpha'_2$, $E'_1$ and $E'_2$ in
$D'$}

\label{hconj-fig07}
\end{center}
\end{figure}

\bca\label{xml-2} $\delta_{D'}(x'_1,x'_2)>
H_1\delta_{D'}(x'_1,x'_3)$.\eca

In this case, we shall exploit the conformal modulus to obtain a homeomorphism. For this, we need some preparation.
Let $\alpha'_1$ be an arc joining $x'_1$ and $x'_2$ in $A'$ (see Figure \ref{hconj-fig07}), and
let $x'_4$ be the first point in $\alpha'_1$ along the direction
from $x'_2$ to $x'_1$ such that (see Figure \ref{hconj-fig07})
\be\label{w-manzi-1}
\delta_{D'}(x'_4,x'_2)= \frac{1}{8}\delta_{D'}(x'_1,x'_2).
\ee
Then we get
\be\label{manzi-4}
\delta_{D'}(x'_4,x'_1)\leq \delta_{D'}(x'_1,x'_2)+\delta_{D'}(x'_4,x'_2)=
\frac{9}{8}\delta_{D'}(x'_1,x'_2).
\ee

In what follows, two claims stated as below are needed.

\bcl\label{manzi-1a} \label{hz-1-1}
$\delta_{D}(x_4,x_2)\geq
H_2\delta_{D}(x_1,x_2)$ with
$\ds H_2=\frac{\psi^{-1}(1/9)}{1+\psi^{-1}(1/9)}$.\ecl

Now, we prove this claim. If $\delta_{D}(x_4,x_2)\geq
\delta_{D}(x_4,x_1)$, then
$$ \delta_{D}(x_4,x_2)\geq \frac{1}{2}\delta_{D}(x_1,x_2),
$$
since $\delta_{D}(x_4,x_2)+\delta_{D}(x_4,x_1)\geq
\delta_{D}(x_1,x_2)$. The claim is true since $H_2<1/2$.

On the other hand, if $\delta_{D}(x_4,x_2)< \delta_{D}(x_4,x_1)$, then by Lemma
\ref{mzz-1} together with (\ref{w-manzi-1}) and (\ref{manzi-4}), we have
$$\frac{1}{9}\leq\frac{\delta_{D'}(x'_4,x'_2)}{\delta_{D'}(x'_4,x'_1)}\leq
\psi\Big(\frac{\delta_{D}(x_4,x_2)}{\delta_{D}(x_4,x_1)}\Big),
$$
whence
$$\delta_{D}(x_4,x_2)\geq\frac{\psi^{-1}(1/9)}{1+\psi^{-1}(1/9)}\delta_{D}(x_1,x_2).
$$
The proof of Claim \ref{manzi-1a} is complete.
\medskip

Let $\alpha'_2$ be an arc
joining $x'_1$ and $x'_3$ in $A'$ (see Figure \ref{hconj-fig07}). Then we have

\bcl\label{ss-1} There exists a point $y'_3$ in $\alpha'_2$ which
satisfies

$(1)$ $\delta_{D'}(x'_1,x'_3)\leq\delta_{D'}(x'_1,y'_3)\leq
2H\delta_{D'}(x'_1,x'_3)$;

$(2)$ $\delta_D(x_1,x_3)\leq \delta_D(x_1,y_3)$; and

$(3)$
 $\alpha'_2[x'_1,y'_3]\subset\overline{\mathbb{B}}(x'_1,2H\delta_{D'}(x'_1,x'_3))$.
\ecl

In order to establish the existence of $y'_3$, we separate the discussions into two parts: $\delta_{D'}(\alpha'_2)\leq 2H\delta_{D'}(x'_1,x'_3)$ and $\delta_{D'}(\alpha'_2)> 2H\delta_{D'}(x'_1,x'_3)$. For the first part, we
let $y'_3=x'_3$. Obviously, $y'_3$ satisfies all requirements in Claim \ref{ss-1}.

For the remaining part, that is, $\delta_{D'}(\alpha'_2)> 2H\delta_{D'}(x'_1,x'_3)$, we see that if
$\alpha'_2[x'_1,x'_3]\cap \mathbb{S}(x'_1,2H\delta_{D'}(x'_1,x'_3))\not=\emptyset$,
we take
$y'_3$ to be the first point in $\alpha'_2$ along the direction from
$x'_1$ to $x'_3$ such that
$$\delta_{D'}(x'_1,y'_3)=2H\delta_{D'}(x'_1,x'_3),
$$
and if $\alpha'_2[x'_1,x'_3]\subset\mathbb{B}(x'_1,2H\delta_{D'}(x'_1,x'_3))$, then let $y'_3\in\alpha'_2$
be such that
$$\delta_{D'}(x'_1,y'_3)>H\delta_{D'}(x'_1,x'_3).
$$
Necessarily, we see that $\delta_D(x_1,x_3)\leq \delta_D(x_1,y_3)$.
Also, the chosen point $y'_3$ satisfies all requirements in Claim \ref{ss-1}. Hence the claim is true.
\medskip

Let us continue the proof of this theorem. Let $E'_1=\alpha'_1[x'_2,x'_4]$ and
$E'_2=\alpha'_2[x'_1,y'_3]$ (see Figure \ref{hconj-fig07}). We need lower
bounds for the quantity $\min\{\diam(E_1), \diam(E_2)\}$ and for the length of
every arc connecting $E'_1$ and $E'_2$ in $D'$, respectively. For this, it follows
from Claim \ref{manzi-1a} that
$$\diam(E_1)\geq \delta_D(x_4,x_2)\geq H_2\delta_D(x_1,x_2)
$$
and from Claim  \ref{ss-1} that
$$\diam(E_2)\geq \delta_D(x_1,y_3)\geq \delta_D(x_1,x_3),
$$
whence the assumption ``$\delta_D(x_1,x_2)>\delta_D(x_1,x_3)$" implies
\beq\label{manzi-8}
\min\{\diam(E_1), \diam(E_2)\}&\geq&
\min\{H_2\delta_D(x_1,x_2), \delta_D(x_1,x_3)\}\\
\nonumber&\geq& H_2\delta_D(x_1,x_3).
\eeq

Since Claim \ref{ss-1}$(3)$ implies $\delta_{D'}(\alpha'_2[x'_1,y'_3])\leq 4H\delta_{D'}(x'_1,x'_3)$,
it follows from the choice of $x_4'$ and the assumption
``$\delta_{D'}(x'_1,x'_2)> H_1\delta_{D'}(x'_1,x'_3)$" that for each
$u'\in E'_1$ and $v'\in E'_2$,
\beq\label{manzi-12}
\delta_{D'}(u',v')&\geq&
\delta_{D'}(x'_1,x'_2)-\delta_{D'}(x'_2,u')-\delta_{D'}(x'_1,v')\\
\nonumber &\geq&
\Big(\frac{7}{8}-\frac{4H}{H_1}\Big)\delta_{D'}(x'_1,x'_2). \eeq

\begin{figure}[htbp]
\begin{center}
\input{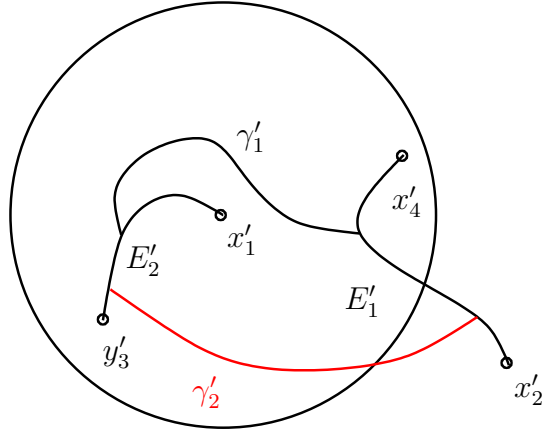}
\caption{The arcs $\gamma'_1\in \Gamma'_1$ and $\gamma'_2\in
\Gamma'_2$ in $D'$}

\label{hconj-fig08}
\end{center}
\end{figure}

The needed lower bounds have been obtained. In order to apply the conformal modulus together with
Theorems \Ref{lemm1--0''} and \Ref{lemm1--xx} to get a homeomorphism, we still need to construct a family of curves in $D'$. Let $\Gamma'$ denote the curves joining $E'_1$ and
$E'_2$ in $D'$. Then (see Figure \ref{hconj-fig08})
$$\Gamma'=\Gamma'_1\cup\Gamma'_2 
,
$$
where
$$\Gamma'_1=\{\gamma'\in\Gamma': \gamma'\subset \mathbb{B}(x'_1, (\delta_{D'}(x'_1,x'_2)\delta_{D'}(x'_1,y'_3))^{\frac{1}{2}})\}
$$
and
$$\Gamma'_2=\{\gamma'\in\Gamma': \gamma'\cap \mathbb{S}(x'_1, (\delta_{D'}(x'_1,x'_2)\delta_{D'}(x'_1,y'_3))^{\frac{1}{2}})\not=\emptyset\}.
$$
Hence we infer from (\ref{manzi-12}), Theorems \Ref{lemm1--0''} and \Ref{lemm1--xx} that
\beq\label{hz-1-7}
{\rm Mod}(E'_1,E'_2;D') &=& {\rm Mod}(\Gamma')\\ \nonumber &\leq& {\rm Mod}(\Gamma'_1)+{\rm Mod}(\Gamma'_2)\\
\nonumber&\leq&\frac{\omega_{n-1}(\delta_{D'}(x'_1,x'_2)\delta_{D'}(x'_1,y'_3))^{\frac{n}{2}}}
{((7/8-4H/H_1)\delta_{D'}(x'_1,x'_2))^n}\\
\nonumber&&+\omega_{n-1}
\Big(\frac{1}{2}\log \frac{\delta_{D'}(x'_1,x'_2)}{\delta_{D'}(x'_1,y'_3)}\Big)^{1-n}\\
\nonumber&\leq& 2^n\omega_{n-1}\Big(\log
\frac{\delta_{D'}(x'_1,x'_2)}{H^2\delta_{D'}(x'_1,x'_3)}\Big)^{1-n},
\eeq
since the statement $(3)$ in Claim \ref{ss-1} shows that for each $\gamma'\in \Gamma'_2$, $\gamma'\cap
\mathbb{S}(x'_1,2H\delta_{D'}(x'_1, x'_3))\not=\emptyset.$ Now, the quasiconformal invariance property of moduli of the families of curves together with
\eqref{manzi-8} guarantees
\begin{eqnarray*}
{\rm Mod}(\Gamma')&\geq& \frac{1}{K}{\rm Mod}(\Gamma)\geq
\frac{1}{K}\varphi\Big(\frac{\delta_D(E_1,E_2)}{\min\{\diam(E_1), \diam(E_2)\}}\Big)\\
\nonumber
&\geq&\frac{1}{K}\varphi\Big(\frac{\delta_D(x_1,x_2)}{H_2\delta_D(x_1,x_3)}\Big),
\end{eqnarray*}
which, combining with (\ref{hz-1-7}), shows that
$$\frac{\delta_{D'}(x'_1,x'_2)}{\delta_{D'}(x'_1,x'_3)}\leq
H^2\exp\bigg({2\Big(\frac{2K\omega_{n-1}}{\tau}\Big)^{\frac{1}{n-1}}}\bigg),
$$
where
$\tau=\varphi\Big(\frac{\delta_D(x_1,x_2)}{H_2\delta_D(x_1,x_3)}\Big)$.
Now, we let
\beq\label{fun-4}
\eta_3(t)=H^2\exp\bigg({2\Big(\frac{2K\omega_{n-1}}{\varphi(t/H_2)}\Big)^{\frac{1}{n-1}}}\bigg)
\eeq
for $t> 1$.

Now, we are ready to conclude the existence of the homeomorphism
$\eta$.
We see from \eqref{sun-1}, \eqref{fun-2} and \eqref{fun-4} that the
homeomorphism
$$\eta(t)=\begin{cases}
\displaystyle  \ds H_1\mu_{11}\eta_1(t) &
\mbox{ if}\,\; \ds t\in(0,1],\\
\displaystyle \ds \max\{\mu_{11}\eta_2(t),36\psi(1)\eta_3(t) \}
 & \mbox{ if}\,\; \ds t\in(1,\infty),
\end{cases}
$$
is the desired one, where
$\mu_{11}=\exp\bigg({2\big(\frac{2K\omega_{n-1}}{\varphi(1/H_2)}\big)^{\frac{1}{n-1}}}\bigg)$.
The proof of Theorem \ref{thm1} is complete. \qed

\section{The equivalence of John domains}\label{sec-4}

In this section, we shall prove Theorem \ref{thm2}. Throughout this
section, we always assume that $D$ and $D'$ are bounded domains in
$\IR^n$, that $f:\, D\to D'$ is a $K$-quasiconformal mapping, that
$D$ is a simply connected $c$-uniform domain, and that
$f:\,\overline{D}\to \overline{D'}$ is continuous. Before the proof of Theorem
\ref{thm2}, we
prove several necessary lemmas. Our first lemma is as follows.

\begin{lem} \label{lem5-A-11a}
Let $x\in D$. Then for each $I(x)\in \Phi(x)$, we have
$$\max\{d_{D'}(x'),\;\dist(x',I'(x))\}\leq \mu_{12}\diam(I'(x)),
$$
where $$\ds
\mu_{12}=4\psi_n^{-1}\left(\frac{1}{K\mu_4}\phi_n\left(10c\Big(1-\exp\big({-\big(\frac{1}{\mu_2}\log\frac{3}{2}\big)
^{\mu_2}}\big)\Big)\right)\right),$$ $\phi_n$ and
$\psi_n$ are from Theorem {\rm \Ref{ThmG}}, $\mu_2$ and $\mu_4$ are from Theorem  {\rm \Ref{ThmF}}
and Theorem  {\rm \Ref{lem0-0}}, respectively. We recall here that $I'(x)=f(I(x))$.
\end{lem}
\bpf
We shall apply the conformal modulus to prove this lemma.

Let $B'_1=\mathbb{B}(x', \frac{1}{2}d_{D'}(x'))$. We need a lower bound for $\diam(B_1)$.
For each
$z'_3\in \partial B'_1$, by the inequality (\ref{eq-2-9'}), we have
$$k_{D'}(x',z'_3)\geq \log\Big(1+\frac{|x'-z'_3|}{d_{D'}(x')}\Big)= \log\frac{3}{2},
$$
which, together with Theorem \Ref{ThmF}, implies that
\be\label{manzi-1}
\max\{k_{D}(x,z_3),(k_{D}(x,z_3))^{\frac{1}{\mu_2}}\}\geq
\frac{1}{\mu_2}k_{D'}(x',z'_3)\geq \frac{1}{\mu_2}\log\frac{3}{2}.
\ee
If $|x-z_3|< \frac{1}{2}d_D(x)$, then $k_D(x, z_3)\leq 1$, and so (\ref{manzi-1}) implies
$$k_{D}(x,z_3)\geq \Big(\frac{1}{\mu_2}\log\frac{3}{2}\Big )^{\mu_2},
$$
whence the inequality (\ref{upperbdk}) leads to
$$\Big(\frac{1}{\mu_2}\log\frac{3}{2}\Big)^{\mu_2}\leq k_D(x,z_3)\leq\log\Big(1+\frac{|z_3-x|}{d_D(x)-|z_3-x|}\Big).
$$
Hence we have proved that for $z_3\in \partial B_1$,
\beq\label{correction-1}
\nonumber
|x-z_3|&\geq &\min\left\{1-\exp\Big({-\big(\frac{1}{\mu_2}\log\frac{3}{2}\big)^{\mu_2}}\Big),\frac{1}{2}\right\}d_{D}(x)\\
\nonumber &=& \left(1-\exp\Big({-\big(\frac{1}{\mu_2}\log\frac{3}{2}\big)^{\mu_2}}\Big)\right)d_{D}(x),
\eeq
which implies
\beq \label{imm-1}
\diam(B_1)\geq 2\left(1-\exp\Big({-\big(\frac{1}{\mu_2}\log\frac{3}{2}\big)^{\mu_2}}\Big)\right)d_{D}(x).
\eeq

Since $I(x)\subset \partial D$, we need to find its replacement in $D$. We construct an arc as follows.
Let $I_1(x)$ be an arc in $D$ which satisfies the following three requirements:
\beq\label{correction-2}
\diam(I(x))\leq 2\diam(I_1(x)),\eeq for each $z\in I_1(x)$,
\beq\label{correction-3}\dist(z,I(x))<\frac{1}{4}\min\{\dist(B_1,I(x)),\diam(I(x))\}
\eeq
and
\beq\label{correction-4}
\dist(z',I'(x))<\frac{1}{4}\min\{\dist(B'_1,I'(x)),\diam(I'(x))\}.
\eeq

In order to apply the conformal modulus to the sets $B_1$ and $I_1(x)$, some preparation is still necessary.
Let $u'_0\in I'_1(x)$ be such that $$\dist(u'_0,B'_1)=\dist(B'_1,I'_1(x)).$$ Then (\ref{correction-4}) implies
\beq\label{correction-14}\dist(B'_1,I'_1(x))\geq\dist(B'_1,I'(x))-\dist(u'_0,I'(x))\geq \frac{3}{4}\dist(B'_1,I'(x)).\eeq

For a comparison between $\diam(I'(x))$ and $\diam(I'_1(x))$, we let $u'_1$ and $u'_2$
in $I'_1(x)$ to be such that $$|u'_1-u'_2|=\diam(I'_1(x)).$$ Then by (\ref{correction-4}), we get
\beq\label{imm-2}
\diam(I'_1(x))&\leq&\dist(u'_1,I'(x))+\dist(u'_2,I'(x))+\diam(I'(x))
\\ \nonumber&\leq& \frac{3}{2}\diam(I'(x)).
\eeq

Moreover, we need two estimates concerning $B_1$ and $I_1(x)$.
For the first estimate, let $z\in I_1(x)$. Then by (\ref{correction-3}) and the
obvious fact ``$\dist(B_1,I(x))\leq 8cd_D(x)$", we have
\beq
\nonumber |z-x|\leq8cd_D(x)+ \dist(z,I(x))\leq 8cd_D(x)+\dist(B_1,I(x))\leq 10cd_D(x),
\eeq
which shows
\beq\label{correction-5}
\dist(I_1(x), B_1)\leq 10cd_D(x),
\eeq
and we see from (\ref{imm-1}) and (\ref{correction-2}) that
\beq\label{correction-1-1} \hspace{.8cm}
\min\{\diam(B_1), \diam(I_1(x))\} &\geq& \min\{\diam(B_1),
\frac{1}{2}\diam(I(x))\}
\\ \nonumber &\geq&
\min\Big \{2-2\exp\Big({-\big(\frac{1}{\mu_2}\log\frac{3}{2}\big)^{\mu_2}}\Big),
\frac{1}{2}\Big\}d_D(x)
\\ \nonumber &\geq&
\left(1-\exp\big({-\big(\frac{1}{\mu_2} \log\frac{3}{2}
\big)^{\mu_2}}\big)\right)d_D(x),\eeq
whence, we obtain
$$B_1\cap I_1(x)=\emptyset.
$$

Now, the quasiconformal invariance property of moduli of the families of curves
implies that
$${\rm Mod}\,(B_1,I_1(x);D)\leq K{\rm Mod}\,(B'_1,I'_1(x);D'),
$$
and so it follows from (\ref{correction-5}), (\ref{correction-1-1}), Theorems \Ref{ThmG}
and \Ref{lem0-0} that
\begin{eqnarray*}
\phi_n\left(10c\Big(1-\exp\big({-\big(\frac{1}{\mu_2}\log\frac{3}{2}\big)^{\mu_2}}\big)\Big)^{-1}\right)
&\leq& \phi_n\Big(\frac{\dist(B_1,I_1(x))}{\min\{\diam(B_1),
\diam(I_1(x))\}}\Big)
\\ \nonumber&\leq& {\rm Mod}\,(B_1,I_1(x);\IR^n)\\
\nonumber&\leq&
 \mu_4{\rm Mod}\,(B_1,I_1(x);D)\\
\nonumber&\leq& K\mu_4{\rm Mod}\,(B'_1,I'_1(x);D')\\
\nonumber&\leq&
K\mu_4{\rm Mod}\,(B'_1,I'_1(x);\IR^n)\\
\nonumber&\leq&
K\mu_4\psi_n\Big(\frac{\dist(B'_1,I'_1(x))}{\min\{\diam(B'_1),
\diam(I'_1(x))\}}\Big),
\end{eqnarray*}
which together with (\ref{correction-14}) and \eqref{imm-2} leads to
\begin{eqnarray*}
\max\{d_{D'}(x'),\dist(x',I'(x))\}& = &\dist(x',I'(x))
\leq  2\dist(B'_1,I'(x))\leq \frac{8}{3}\dist(B'_1,I'_1(x))  \\ &\leq &
\frac{8}{3}\psi_n^{-1}\left(\frac{1}{K\mu_4}\phi_n\left(10c\Big(1-\exp\big({-\big(\frac{1}{\mu_2}\log\frac{3}{2}\big)^{\mu_2}}\big)\Big)\right)\right)\\
&& \times\min\{\diam(B'_1),
\diam(I'_1(x))\}
\\ &\leq &
\frac{8}{3}\psi_n^{-1}\left(\frac{1}{K\mu_4}\phi_n\left(10c\Big(1-\exp\big({-\big(\frac{1}{\mu_2}\log\frac{3}{2}\big)^{\mu_2}}\big)\Big)\right)\right)\\
&& \times\diam(I'_1(x))\\ &\leq &
4\psi_n^{-1}\left(\frac{1}{K\mu_4}\phi_n\left(10c\Big(1-\exp\big({-\big(\frac{1}{\mu_2}\log\frac{3}{2}\big)^{\mu_2}}\big)\Big)\right)\right)\\
&& \times\diam(I'(x)).
\end{eqnarray*}
Hence the proof of the lemma is complete. \epf

In Lemma \ref{lem5-A-11a}, a lower bound for the ratio $\frac{\diam(I'(x))}{d_{D'}(x')}$ has been proved.
In the following lemmas, we shall obtain some other upper bounds for $\frac{\diam(I'(x))}{d_{D'}(x')}$ under
different conditions. All these bounds are needed later on.

\begin{lem} \label{lem5-A-11}
Suppose there are constants $b\geq 1$ and $\alpha\leq 1$ such that
$$\frac{\diam(P')}{\diam(Q')}\leq b\Big(\frac{\diam(P)}{\diam(Q)}\Big)^{\alpha}
$$
for all components $P\subset Q$, where $P\in \Phi(u)$ and $Q\in \Phi(v)$ for $u$, $v\in D$ (here the case $u=v$ inclusive).
Then for $x\in D$ and each $I(x)\in \Phi(x)$, we have
$$\diam(I'(x))\leq \mu_{13}d_{D'}(x'),
$$
where
$$\ds \mu_{13}=\max\left \{72(\exp(\mu_2)-1),
24(\exp(\mu_2)-1)\psi_n^{-1}\Big(\frac{\phi_n(160c^2(4b
)^{\frac{1}{\alpha}})}{K\mu_4}\Big)\right \}.
$$
\end{lem}

\bpf Suppose on the contrary that there exist some $x\in D$ and
$I(x)\in \Phi(x)$ such that
\be\label{manzi-3-10}
\diam(I'(x))> \mu_{13}d_{D'}(x').
\ee

Under this assumption, we shall exploit the conformal modulus to get a contradiction.
At first, we do some preparation.

Obviously, the assumption \eqref{manzi-3-10} assures that there must exist a continuum $P_0\subset
I(x)$ such that
\be\label{ym-1}
\diam(P'_0)> \frac{1}{4}\diam(I'(x))\;\; \mbox{and}\;\; \dist(x', P'_0)\geq
\frac{1}{6}\diam(I'(x)).
\ee
Then there must exist some point $u\in D$ such that
\be\label{ymws-1}
\diam(P_0)\leq 3cd_D(u)\;\; \mbox{and}\;\; d_D(u)\leq\dist(u, P_0)\leq \frac{1}{2}\diam(P_0),
\ee and thus
$$P_0\subset \mathbb{B}(u,8cd_D(u)),
$$
whence \eqref{ymws-1} guarantees that there exists some $P\in \Phi(u)$ such that
$$P_0\subset P.
$$
Hence the assumption in the lemma implies
$$\frac{1}{4}<\frac{\diam(P'_0)}{\diam(I'(x))}\leq \frac{\diam(P')}{\diam(I'(x))}\leq
b\Big(\frac{\diam(P)}{\diam(I(x))}\Big)^{\alpha},
$$
which shows that
\be\label{ym-2}
8c\diam(P_0)\geq \diam(P)> \Big(\frac{1}{4b}\Big)^{\frac{1}{\alpha}}\diam(I(x)),
\ee
since (\ref{ymws-1}) implies $\diam(P)\leq 16cd_D(u)\leq 8c\diam(P_0)$.

Let $B_2=\mathbb{B}(x, \frac{1}{2}d_D(x))$. We need estimates on $\diam(B'_2)$ and
$\dist(B'_2,P'_0)$. For the first estimate, we let $y_3\in
\partial B_2$. Apparently, $k_D(x,y_3)\leq 1,$ and so Theorem
\Ref{ThmF} and (\ref{eq-2-9'}) lead to
$$\log\Big(1+\frac{|x'-y'_3|}{d_{D'}(x')}\Big)\leq k_{D'}(x',y'_3)\leq
\mu_2(k_D(x,y_3))^{\frac{1}{\mu_2}}\leq \mu_2,
$$
whence
\be\label{ym-3}
\diam(B'_2)\leq 2(\exp({\mu_2})-1)d_{D'}(x'),
\ee
and for the second one, we see from (\ref{manzi-3-10}), (\ref{ym-1}) and (\ref{ym-3}) that
\beq\label{manzi-3-12}
\dist(B'_2,P'_0)&\geq& \dist(x',P'_0)-\diam(B'_2)\\ \nonumber&\geq&
\Big(\frac{\mu_{13}}{6}-2(\exp({\mu_2})-1)\Big)d_{D'}(x').
\eeq

Since $P_0\subset \partial D$, we need its replacement of $P_0$ in $D$. We construct an arc as follows:
Let $P_1\subset D$ be an arc which satisfies the following three requirements:
\beq\label{correction-7}
\diam(P_0)\leq 2\diam(P_1),
\eeq
for each $z\in P_1$,
\beq\label{correction-8}
\dist(z,P_0)<\frac{1}{4}\min\{\dist(B_2,P_0),\diam(P_0)\}
\eeq
and
\beq\label{correction-9}
\dist(z',P'_0)<\frac{1}{4}\min\{\dist(B'_2,P'_0),\diam(P'_0)\}.
\eeq
Then we establish several inequalities related to $P_1$ and its image $P'_1$. For the first
inequality, we let $u'_3\in P'_1$ be such that $\dist(u'_3,B'_2)=\dist(B'_2,P'_1)$.
Then by (\ref{manzi-3-12}) and (\ref{correction-9}), we have
\beq\label{correction-13}
\dist(B'_2,P'_1)&\geq& \dist(B'_2,P'_0)-\dist(u'_3,P'_0)\geq
\frac{3}{4}\dist(B'_2,P'_0)\\ \nonumber&\geq& \Big(\frac{\mu_{13}}{8}-\frac{3}{2}(\exp({\mu_2})-1)\Big)d_{D'}(x').
\eeq

For the second one, we let $u_4\in P_1$ be such that $\dist(x,P_1)=|x-u_4|$ and
$u_5\in P_0$ such that $\dist(u_4,P_0)=|u_4-u_5|$. Then by (\ref{correction-8}), we see that
\beq\label{correction-17}
\dist(B_2,P_1)&\leq& \dist(x,P_1)\\
\nonumber&\leq&|x-u_5|+|u_4-u_5|\\
\nonumber&\leq&8cd_D(x)+\frac{1}{4}\dist(B_2,P_0)\\
\nonumber&\leq&8cd_D(x)+\frac{1}{4}|x-u_5|\\
\nonumber&\leq& 10cd_D(x),
\eeq
since $P_0\subset I(x) \subset\mathbb{B}(x,8cd_D(x))$, i.e., $|x-u_5|\leq 8cd_D(x)$.

For the third inequality, we see from (\ref{ym-2}) and (\ref{correction-7}) that
\beq\label{correction-11}
\diam(P_1)&\geq&\frac{1}{2}\diam(P_0)\geq \frac{1}{16c}\cdot\Big(\frac{1}{4b}\Big)^{\frac{1}{\alpha}}\diam(I(x))
\\
\nonumber&\geq& \frac{1}{16c}\cdot\Big(\frac{1}{4b}\Big)^{\frac{1}{\alpha}}d_D(x).
\eeq

Obviously,
$$B_2\cap P_1=\emptyset.
$$
Now, the quasiconformal invariance property of moduli of the families of curves implies
$${\rm Mod}\,(B_2,P_1;D)\leq K{\rm Mod}\,(B'_2,P'_1;D'),
$$
whence it follows from Theorems \Ref{ThmG} and \Ref{lem0-0} together with (\ref{ym-3}),
(\ref{correction-13}), (\ref{correction-17}) and (\ref{correction-11}) that
\begin{eqnarray*}
\phi_n\Big(160c^2 (4b )^{\frac{1}{\alpha}}\Big)
&\leq& \phi_n\Big(\frac{10cd_D(x)}{\min\{\diam(B_2),
\diam(P_1)\}}\Big)\\
&\leq & \phi_n\Big(\frac{\dist(B_2,P_1)}{\min\{\diam(B_2),
\diam(P_1)\}}\Big)
\\ &\leq& {\rm Mod}\,(B_2,P_1;\IR^n)
\leq \mu_4{\rm Mod}\,(B_2,P_1;D)\\
&\leq& K\mu_4{\rm Mod}\,(B'_2,P'_1;D')\\
&\leq& K\mu_4\psi_n\Big(\frac{\dist(B'_2,P'_1)}{\min\{\diam(B'_2),
\diam(P'_1)\}}\Big)\\
&\leq& K\mu_4\psi_n\Big(
\frac{\frac{1}{8}\mu_{13}-\frac{3}{2}(\exp({\mu_2})-1)}{2(\exp({\mu_2})-1)}\Big)\\
&<& K\mu_4\psi_n\Big(\frac{\mu_{13}}{24(\exp({\mu_2})-1)}\Big)\\
&<& \phi_n\Big(160c^2 (4b )^{\frac{1}{\alpha}}\Big).
\end{eqnarray*}
 This obvious contradiction completes the proof. \epf

\begin{figure}[htbp]
\begin{center}
\input{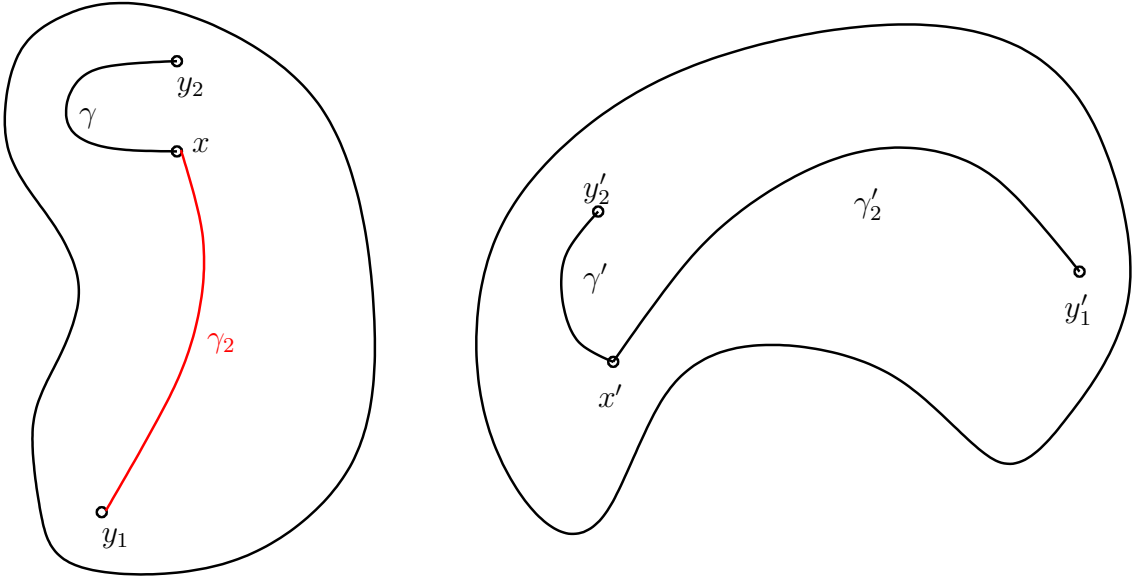}
\caption{The arcs $\gamma$, $\gamma_2$ in $D$, and their images
under $f$ in $D'$}

\label{hconj-fig09}
\end{center}
\end{figure}

\begin{lem}\label{Sun-1}
Suppose $f:\,(D,\delta_D)\to (D', \delta_{D'})$ is
$\eta$-quasisymmetric. Then for $x\in D$ and each $I(x)\in \Phi(x)$,
we have
$$\diam\big(I'(x)\big)\leq \mu_{14}d_{D'}(x'),
$$
where $\ds \mu_{14}=3\eta(16c^2)(\exp(\mu_2)-1)$.
\end{lem}
\bpf
For $x\in D$, let $y'_1\in  D'$ be such that (see Figure \ref{hconj-fig09})
\be\label{manzi-3-1}
y_1\in \mathbb{B}(x,8cd_D(x))\;\; \mbox{and}\;\;|x'-y'_1|\geq \frac{1}{3}\diam(I'(x)),
\ee
and let $y_2\in D$ be such that (see Figure \ref{hconj-fig09})
$$|y_2-x|=\frac{1}{2}d_D(x).
$$
Obviously, $k_D(x,y_2)\leq 1,$ and so Theorem \Ref{ThmF} and (\ref{eq-2-9}) imply that
$$\log\Big(1+\frac{\ell(\gamma')}{d_{D'}(x')}\Big)\leq
k_{D'}(x',y'_2)\leq\mu_2\max\{k_D(x,y_2),(k_D(x,y_2))^{\frac{1}{\mu_2}}\}
\leq \mu_2,
$$
where $\gamma'$ is a quasihyperbolic geodesic joining $x'$ and $y'_2$ in $D'$ (see Figure \ref{hconj-fig09}), and then
\be\label{mlx-3-1}
\ell(\gamma')\leq (\exp({\mu_2})-1)d_{D'}(x').
\ee

With the aid of $\ell(\gamma')$, we shall finish the proof. For this, we let $\gamma_2$ be a double
$c$-cone arc joining $y_1$ and $x$ in $D$ (see Figure \ref{hconj-fig09}). The
existence of $\gamma_2$ follows from the assumption that ``$D$ is
$c$-uniform". Then
\be\label{mlx-3-2}
\delta_{D}(x,y_1)\leq \ell(\gamma_2)\leq c|y_1-x|,
\ee
and thus the hypotheses on $f$ in the lemma, (\ref{manzi-3-1}) and (\ref{mlx-3-2}) show that
$$\frac{\diam(I'(x))}{\delta_{D'}(x',y'_2)}\leq \frac{3\delta_{D'}(x',y'_1)}{\delta_{D'}(x',y'_2)}
\leq 3\eta\Big(\frac{\delta_{D}(x,y_1)}{\delta_{D}(x,y_2)}\Big)\leq
3\eta(16c^2),
$$
which, together with (\ref{mlx-3-1}), leads to
$$\diam(I'(x))\leq 3\eta(16c^2)\delta_{D'}(x',y'_2)\leq 3\eta(16c^2)\ell(\gamma')
\leq 3\eta(16c^2)(\exp({\mu_2})-1)d_{D'}(x'),
$$
from which the proof follows by taking
$\mu_{14}=3\eta(16c^2)(\exp({\mu_2})-1)$.\epf

\begin{lem} \label{lem5-A-11-1}
Suppose there are constants $\epsilon_0\geq 1$ and $\alpha\leq 1$
such that for $x, w \in D$ with $|x-w|\leq 8cd_D(x)$ and
$d_D(w)\leq 2cd_D(x)$, the inequality
$$a_f(w)\leq \epsilon_0 a_f(x)\Big(\frac{d(x)}{d(w)}\Big)^{1-\alpha}
$$
holds. Then for $I(x)\in \Phi(x)$, we have
$$\diam(I'(x))\leq \mu_{16} d_{D'}(x'),
$$
where
 $$\ds \mu_{16}=3(1+4c\mu_{15})^{\mu_1\mu_2}\bigg(1+(\exp({\mu_2})-1)
\bigg(\frac{20c}{5c-1}+\frac{2}{5c-1}\psi_n^{-1}\Big(\frac{\phi_n(c\mu_{15})}{K\mu_4}\Big)\bigg)\bigg),
$$ $\ds \mu_{15}=(5c\epsilon_0\mu_5^2)^{\frac{3}{\alpha}},$ $\mu_1$ and $\mu_5$ are from Theorem {\rm \Ref{ThmF'}} and Theorem {\rm \Ref{mzxl-3-1}}, respectively.
\end{lem}

\begin{figure}[htbp]
\begin{center}
\input{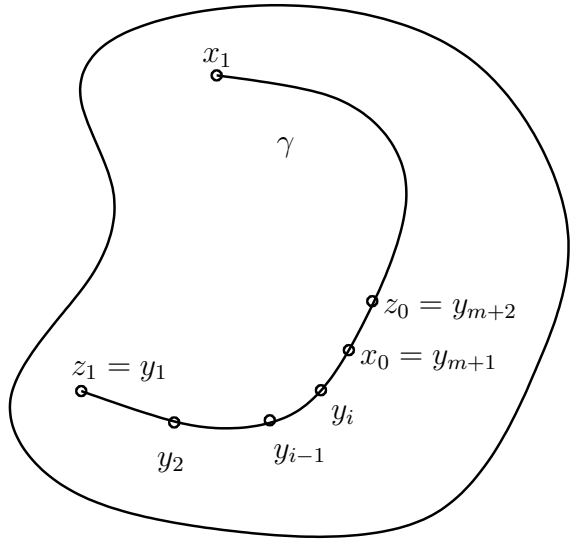}
\caption{The arc $\gamma$ and the related points in $D$
}
\label{hconj-fig10}
\end{center}
\end{figure}

\bpf For $x_1\in D$, we take $z_1\in \mathbb{B}(x_1,8cd_D(x_1))\cap
D$ such that (see Figure \ref{hconj-fig10})
\be\label{xtm-1}
d_D(z_1)\leq
\Big(\frac{1}{5c\epsilon_0\mu_5^2}\Big)^{\frac{3}{\alpha}}d_D(x_1)
\;\;\mbox{and}\;\; |z'_1-x'_1|\geq \frac{1}{3}\diam(I'(x_1)).
\ee
Then we have
\beq\label{wedn-1}
|z_1-x_1|>\frac{1}{2}d_D(x_1),
\eeq
because otherwise,
$$d_D(z_1)\geq d_D(x_1)-|z_1-x_1|\geq \frac{1}{2}d_D(x_1),
$$
which contradicts (\ref{xtm-1}).

Let $\gamma$ denote a double $c$-cone arc joining $z_1$ and $x_1$ in
$D$ (see Figure \ref{hconj-fig10}). Then
\be\label{xtm-2}
\ell(\gamma)\leq c|z_1-x_1|\leq 4cd_D(x_1).
\ee
Let $z_0$ be the midpoint of $\gamma$ with respect to the arc length (see Figure \ref{hconj-fig10}) and
$\mu_{15}=(5c\epsilon_0\mu_5^2)^{\frac{1}{\alpha}}$, where $\mu_5$
is the same constant as in Theorem \Ref{mzxl-3-1}. In the following, we partition $\gamma$. Since by
\eqref{wedn-1},
$$d_D(z_0)\geq \frac{1}{2c}|z_1-x_1|>\frac{1}{4c}d_D(x_1),
$$
it follows from \eqref{xtm-1} that there must exist an
integer $m$~~$(m\geq 2)$ such that
\be\label{xtmz-3-1}
\mu_{15}^{m}\, d_D(z_1) \leq d_D(z_0)< \mu_{15}^{m+1}\, d_D(z_1).
\ee
We use $x_0$ to denote the first point in $\gamma$ from $z_1$ to $z_0$ satisfying (see Figure \ref{hconj-fig10})
$$d_D(x_0)=\mu_{15}^{m}\, d_D(z_1).
$$
Let $y_1=z_1$. Then we choose points $y_2, \ldots, y_{m+1}$ in $\gamma[z_1,z_0]$
such that for each $i\in \{2,\ldots,m+1\}$, $y_i$ is the first point
from $z_1$ to $z_0$ with
\be\label{hws-eq(4.3)}
d_D(y_i)=\mu_{15}^{i-1}\, d_D(y_1).
\ee

Obviously, $y_{m+1}=x_0$. If $x_0\not= z_0$, then we let $y_{m+2}=z_0$, and thus we get a partition of
$\gamma$ (see Figure \ref{hconj-fig10}).
It easily follows from
\eqref{xtmz-3-1} that
\beq \label{bonn-4}
d_D(y_{m+1})> \frac{1}{\mu_{15}}d_D(z_0)\geq \frac{1}{2c\mu_{15}}\ell(\gamma).
\eeq
We need estimates on $d_{D'}(y'_{m+1})$ and $|y'_{m+1}-x'_1|$ in terms of $d_{D'}(x'_1)$, as well as
an upper bound for $\frac{d_{D'}(y'_{i-1})}{d_{D'}(y'_i)}$.
It follows from Theorem \Ref{ThmF'}, (\ref{xtm-2}) and \eqref{bonn-4} that
$$k_D(y_{m+1},x_1)\leq \mu_1\log\Big(1+\frac{|y_{m+1}-x_1|}{\min\{d_D(y_{m+1}),d_D(x_1)\}}\Big)
\leq\mu_1\log(1+4c\mu_{15}),
$$
and then by Theorem \Ref{ThmF} and \eqref{eq-2-9'}, we have

\vspace{6pt}
$\ds \max\left \{\log\Big(1+\frac{|y'_{m+1}-x'_1|}{d_{D'}(x'_1)}\Big),\log \frac{d_{D'}(y'_{m+1})}{d_{D'}(x'_1)}
\right \}$
\begin{eqnarray*}
&\leq&
k_{D'}(y'_{m+1},x'_1)\\
\nonumber &\leq&\mu_2\max\{k_{D}(y_{m+1},x_1),(k_{D}(y_{m+1},x_1))^{\frac{1}{\mu_2}}\}\\
\nonumber &\leq&\mu_1\mu_2\log(1+4c\mu_{15}),
\end{eqnarray*}
whence
\be\label{xtm-8}
d_{D'}(y'_{m+1})\leq (1+4c\mu_{15})^{\mu_1\mu_2}d_{D'}(x'_1)
\ee
and
\be\label{xtm-9}
|y'_{m+1}-x'_1|\leq \big((1+4c\mu_{15})^{\mu_1\mu_2}-1\big)d_{D'}(x'_1).
\ee
Since for each $i\in \{2,\ldots,m +1\}$,
$$|y_{i-1}-y_i|\leq cd_D(y_i)\;\; \mbox{and}\;\; d_D(y_{i-1})\leq d_D(y_i),
$$
we infer from the assumptions in the lemma that
$$a_f(y_{i-1})\leq \epsilon_0 a_f(y_i)\Big(\frac{d_D(y_i)}{d_D(y_{i-1})}\Big)^{1-\alpha},
$$
which, together with Theorem \Ref{mzxl-3-1} and \eqref{hws-eq(4.3)}, shows that
\beq\label{xtm-3}
d_{D'}(y'_{i-1})&\leq& \mu_5\,a_f(y_{i-1})d_D(y_{i-1})\\ \nonumber&\leq& \mu_5 \epsilon_0
a_f(y_i)\Big(\frac{d_D(y_i)}{d_D(y_{i-1})}\Big)^{1-\alpha}d_D(y_{i-1})\\ \nonumber&\leq& \mu_5^2 \epsilon_0
a_f(y_i)\Big(\frac{d_D(y_i)}{d_D(y_{i-1})}\Big)^{-\alpha}d_{D'}(y'_i)
 \\ \nonumber&\leq&\frac{1}{5c}d_{D'}(y'_i).
\eeq

To complete the proof, we still need to prove the following claim.

\bcl\label{xtmz-1} For each $i\in
\{2,\ldots,m+1\}$, we have
$$|y'_i-y'_{i-1}|\leq 2(\exp({\mu_2})-1)\bigg(2+\frac{1}{5c}\psi_n^{-1}
\Big(\frac{\phi_n(c\mu_{15})}{K\mu_4}\Big)\bigg)d_{D'}(y'_i).
$$
\ecl

We shall get a proof of this lemma by using the conformal modulus. For this, we need to construct
a finite sequence of disjoint balls in $D$. For each
$i\in \{2,\ldots,m+1\}$, we let
$$B_i=\mathbb{B}(y_i, \frac{1}{2}d_{D}(y_i)).
$$
Then for $w\in \partial B_i$, $k_D(y_i,w)\leq 1,$ and so
Theorem \Ref{ThmF} and (\ref{eq-2-9'}) imply that
$$\log\Big(1+\frac{|y'_i-w'|}{d_{D'}(y'_i)}\Big)\leq k_{D'}(y'_i,w')\leq
\mu_2(k_D(y_i,w))^{\frac{1}{\mu_2}}\leq \mu_2,
$$
whence
\be\label{xtm-4}
\diam(B'_i)\leq 2(\exp({\mu_2})-1)d_{D'}(y'_i).
\ee
It easily follows from \eqref{hws-eq(4.3)} and the inequality
$$|y_i-y_{i-1}|\geq d_D(y_i)-d_D(y_{i-1}) =\Big(1-\frac{1}{\mu_{15}}\Big)d_D(y_i)
$$
that for each $i\in \{2,\ldots,m\}$,
$$B_i\cap B_{i-1}=\emptyset.
$$
Then the quasiconformal invariance property of moduli of the families of curves
implies that for each $i\in \{2,\ldots,m\}$,
$${\rm Mod}\,(B_{i-1},B_i;D)\leq K{\rm Mod}\,(B'_{i-1},B'_i;D'),
$$
and so we deduce from Theorems \Ref{ThmG} and \Ref{lem0-0} together with
(\ref{hws-eq(4.3)}) that
\begin{eqnarray*}
\phi_n(c\mu_{15})
&\leq& \phi_n\Big(\frac{cd_D(y_i)}{\min\{\diam(B_{i-1}),
\diam(B_i)\}}\Big)\\
&\leq &\phi_n\Big(\frac{\dist(B_{i-1},B_i)}{\min\{\diam(B_{i-1}),\diam(B_i) \}}\Big)\\
\nonumber&\leq& {\rm Mod}\,(B_{i-1},B_i;\IR^n)\leq
 \mu_4{\rm Mod}\,(B_{i-1},B_i;D)\leq K\mu_4{\rm Mod}\,(B'_{i-1},B'_i;D')\\
\nonumber&\leq&
K\mu_4\psi_n\Big(\frac{\dist(B'_{i-1},B'_i)}{\min\{\diam(B'_{i-1}),
\diam(B'_i)\}}\Big),\\
\end{eqnarray*}
which, together with (\ref{xtm-3}) and (\ref{xtm-4}), shows that
\begin{eqnarray*}
\dist(B'_{i-1},B'_i) &\leq&\psi_n^{-1}\Big(\frac{\phi_n(c\mu_{15})}{K\mu_4}\Big)\min\{\diam(B'_{i-1}),
\diam(B'_i)\}\\ \nonumber&\leq&
\frac{2(\exp({\mu_2})-1)}{5c}\psi_n^{-1}\Big(\frac{\phi_n(c\mu_{15})}{K\mu_4}\Big)d_{D'}(y'_i),
\end{eqnarray*}
and thus by (\ref{xtm-4}),
\begin{eqnarray*}
|y'_i-y'_{i-1}|&\leq&\dist(B'_{i-1},B'_i)+2\max\{\diam(B'_{i-1}), \diam(B'_i)\} \\ \nonumber&\leq&
\frac{2(\exp({\mu_2})-1)}{5c}\psi_n^{-1}\Big(\frac{\phi_n(c\mu_{15})}{K\mu_4}\Big)d_{D'}(y'_i)+4(\exp({\mu_2})-1)d_{D'}(y'_i)\\
\nonumber&=&
2(\exp({\mu_2})-1)\bigg(2+\frac{1}{5c}\psi_n^{-1}\Big(\frac{\phi_n(c\mu_{15})}{K\mu_4}\Big)\bigg)d_{D'}(y'_i).
\end{eqnarray*}
The proof of Claim \ref{xtmz-1} is complete.\medskip

We are ready to conclude the proof. It follows from Claim
\ref{xtmz-1}, (\ref{xtm-1}), (\ref{xtm-8}), (\ref{xtm-9}) and (\ref{xtm-3}) that
\begin{eqnarray*}
\diam(I'(x))&\leq& 3
|y'_1-x'_1|\leq 3(|y'_1-y'_2|+\dots+|y'_{m+1}-x'_1|)\\
\nonumber
&\leq&\frac{30c(\exp({\mu_2})-1)}{5c-1}\bigg(2+\frac{1}{5c}\psi_n^{-1}\Big(\frac{\phi_n(c\mu_{15})}
{K\mu_4}\Big)\bigg)d_{D'}(y'_{m+1})\\ \nonumber &&
+3(1+4c\mu_{15})^{\mu_1\mu_2}d_{D'}(x'_1)
\\
\nonumber
&\leq&3(1+4c\mu_{15})^{\mu_1\mu_2}\bigg(1+(\exp({\mu_2})-1)
\bigg(\frac{20c}{5c-1}\\ \nonumber && +\frac{2}{5c-1}\psi_n^{-1}
\Big(\frac{\phi_n(c\mu_{15})}{K\mu_4}\Big)\bigg)\bigg)d_{D'}(x'_1).
\end{eqnarray*}
By taking
$$\mu_{16}=3(1+4c\mu_{15})^{\mu_1\mu_2}\bigg(1+(\exp({\mu_2})-1)
\bigg(\frac{20c}{5c-1}+\frac{2}{5c-1}\psi_n^{-1}
\Big(\frac{\phi_n(c\mu_{15})}{K\mu_4}\Big)\bigg)\bigg),
$$
we see that the lemma is true.\epf

\subsection{ The proof of Theorem \ref{thm2}}

\begin{figure}[htbp]
\begin{center}

\input{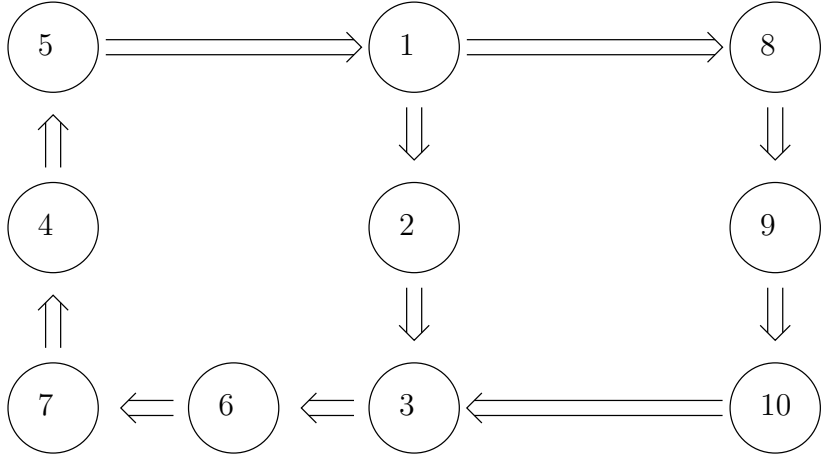}

\caption{The route of the proof of Theorem \ref{thm2}}

\label{hconj-fig23}
\end{center}
\end{figure}

Obviously, it suffices to prove two groups of implications (see Figure \ref{hconj-fig23}):

\noindent (I)\quad \eqref{thm2-3} $\Longrightarrow$ \eqref{thm2-6}
$\Longrightarrow$ \eqref{thm2-7} $\Longrightarrow$ \eqref{thm2-4}
$\Longrightarrow$ \eqref{thm2-5} $\Longrightarrow$ \eqref{thm2-1}
$\Longrightarrow$ \eqref{thm2-2} $\Longrightarrow$ \eqref{thm2-3}; and

\noindent (II)\quad \eqref{thm2-1} $\Longrightarrow$ \eqref{thm2-8} $\Longrightarrow$
\eqref{thm2-9} $\Longrightarrow$ \eqref{thm2-10}  $\Longrightarrow$
\eqref{thm2-3}.

\subsection{ \eqref{thm2-3} $\Longrightarrow$ \eqref{thm2-6}
$\Longrightarrow$ \eqref{thm2-7} $\Longrightarrow$ \eqref{thm2-4}
$\Longrightarrow$ \eqref{thm2-5} $\Longrightarrow$ \eqref{thm2-1}
$\Longrightarrow$ \eqref{thm2-2} $\Longrightarrow$ \eqref{thm2-3}}\
\

The implications \eqref{thm2-1}  $\Longrightarrow$ \eqref{thm2-2}
$\Longrightarrow$ \eqref{thm2-3} easily follow from \cite[Lemma 2.6]{Gem}, \cite[Lemma 3.9]{RJ},
 \cite[Theorem 2.20]{Vai0} and Remark \ref{xt-11}. Lemma \ref{lem5-A-11} shows
that the implication \eqref{thm2-7} $\Longrightarrow$ \eqref{thm2-4}
is true. Therefore, it remains to prove four implications as follows.

\subsubsection{\eqref{thm2-4} $\Longrightarrow$ \eqref{thm2-5}}\label{Sun-2}\
\

 For $w$ and $x$ in $D$ with $|w-x|\leq 8cd_D(x)$, we need to prove
 \beq \label{monkey-1}
 \delta_{D'}(x',w')\leq b_2d_{D'}(x').
 \eeq

\begin{figure}[htbp]
\begin{center}
\small
\input{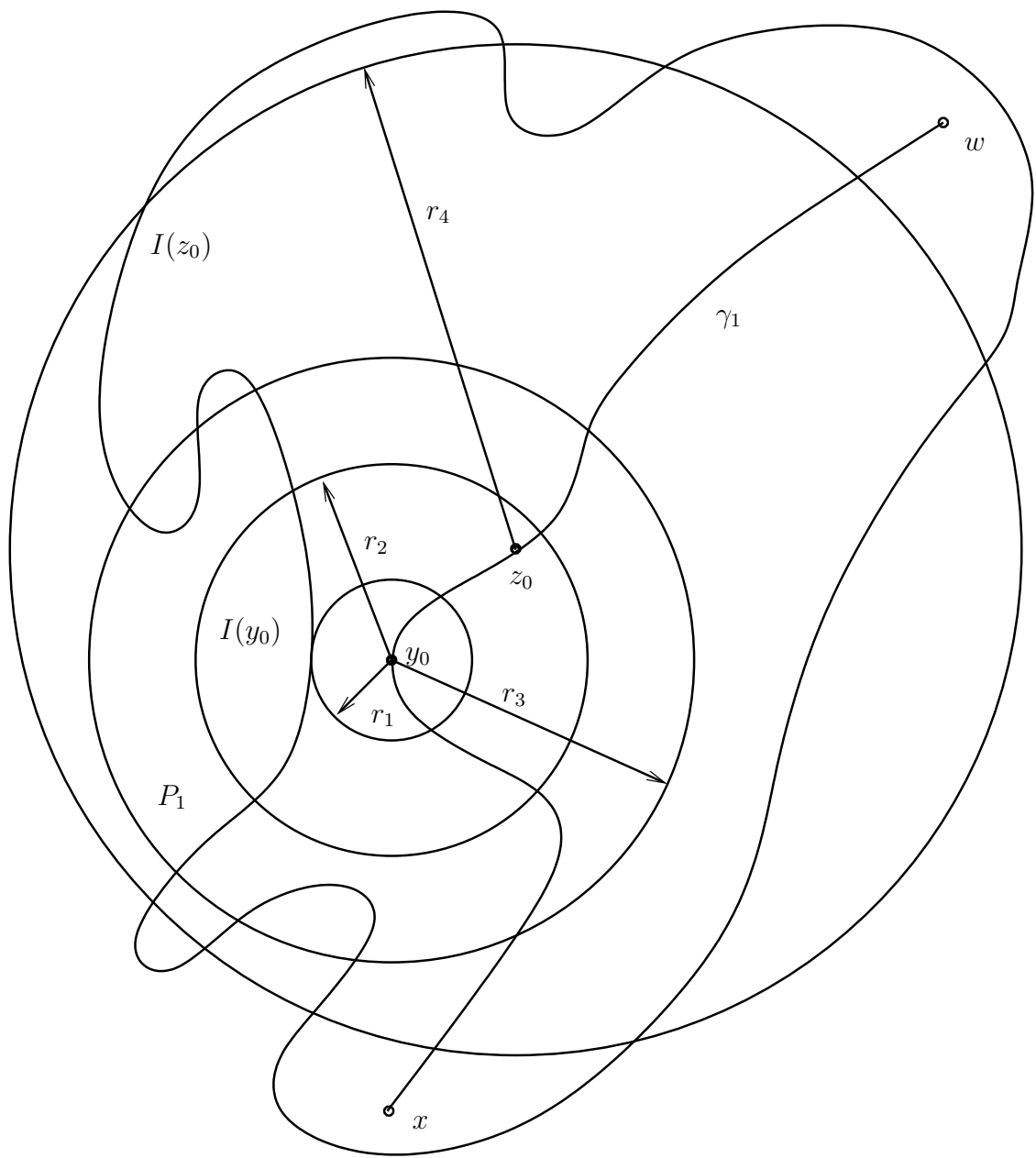}
\normalsize
\caption{The double cone arc $\gamma_1$, the points $x,w,z_0,y_0$, and the continua $I(y_0)\subset P_1\subset I(z_0)$, where $r_1=d_D(y_0)$, $r_2=8cd(y_0)$, $r_3=6cd_D(z_0)$ and $r_4=8cd_D(z_0)$}
\label{hconj-new12b}
\end{center}
\end{figure}

 Without loss of generality, we assume that $x\not=w$.
Let $\gamma_1$ be a double $c$-cone arc joining $w$ and $x$ in $D$,
and let $z_0$ be the midpoint of $\gamma_1$ with respect to the arc
length (see Figure \ref{hconj-new12b}). Obviously,
$$\ell(\gamma_1)\leq c|w-x|\leq8c^2d_D(x).
$$

\noindent Then it follows from Theorem \Ref{ThmF'} that
$$k_D(z_0,x)\leq\mu_1\log\Big(1+\frac{|x-z_0|}{\min\{d_D(z_0),d_D(x)\}}\Big)\leq \mu_1\log(1+4c^2),
$$
whence Theorem \Ref{ThmF} and \eqref{eq-2-9'} lead to
$$\log \frac{d_{D'}(z'_0)}{d_{D'}(x')}\leq
k_{D'}(z'_0,x')\leq\mu_2\max\{k_{D}(z_0,x),(k_{D}(z_0,x))^{\frac{1}{\mu_2}}\} \leq \mu_1\mu_2\log(1+4c^2),
$$
and so
\be\label{whuang-3-4}
d_{D'}(z'_0)\leq (1+4c^2)^{\mu_1\mu_2}d_{D'}(x').
\ee

We choose $y'_0\in\gamma'_1$ such that
\be\label{whuang-3-3}
|y'_0-z'_0|\geq \frac{1}{2}\delta_{D'}(w',x').
\ee
Then
\beq\label{happy-1}
|y_0-z_0|\leq cd_D(z_0).
\eeq
We continue the proof by considering two cases.
If $d_D(y_0)\geq \frac{1}{2}d_D(z_0)$, then by Theorem
\Ref{ThmF'} and (\ref{happy-1}), we have
$$k_D(z_0,y_0)\leq\mu_1\log\Big(1+\frac{|z_0-y_0|}{\min\{d_D(z_0),d_D(y_0)\}}\Big)\leq
\mu_1\log(1+2c),
$$
which, together with \eqref{eq-2-9'} and Theorem \Ref{ThmF}, yields
that
\begin{eqnarray*}
\log \frac{|y'_0-z'_0|}{d_{D'}(z'_0)}
&\leq& k_{D'}(y'_0,z'_0)\leq\mu_2\max\{k_{D}(z_0,y_0),(k_{D}(z_0,y_0))^{\frac{1}{\mu_2}}\}\\
\nonumber &\leq&\mu_1\mu_2\log(1+2c).
\end{eqnarray*}
Hence we infer from (\ref{whuang-3-4}) and \eqref{whuang-3-3} that
\beq\label{Thur-3}
\delta_{D'}(w',x') &\leq&2|y'_0-z'_0|
 \leq 2(1+2c)^{\mu_1\mu_2}d_{D'}(z'_0)\\ \nonumber&\leq&
2(8c^3+4c^2+2c+1)^{\mu_1\mu_2}d_{D'}(x').
\eeq

On the other hand, if $d_D(y_0)< \frac{1}{2}d_D(z_0)$, then it follows from \eqref{happy-1} that
$$\mathbb{B}(y_0,6cd_D(z_0))\subset\mathbb{B}(z_0, 8cd_D(z_0)),
$$
and so there exists some continuum $P_1$ such that (see Figure \ref{hconj-new12b})
\begin{enumerate}
\item $P_1\subset \partial D\cap \mathbb{B}(y_0,
6cd_D(z_0));$
\item $P_1\cap \mathbb{S}(y_0,
d_D(y_0))\not=\emptyset;$
\item $\diam(P_1)\geq 2d_D(z_0);$
\item there exist
$I(z_0)\in \Phi(z_0)$ and $I(y_0)\in
\Phi(y_0)$ such that $I(y_0)\subset
P_1\subset I(z_0).$
\end{enumerate}

The combination of (\ref{whuang-3-4}), Lemma \ref{lem5-A-11a} and
the condition \eqref{thm2-4} in the theorem leads to
\beq\nonumber
|y'_0-z'_0|&\leq& \dist(y'_0,
I'(y_0))+\dist(z'_0, I'(z_0))+\diam(I'(z_0))\\
\nonumber &\leq& (2\mu_{12}+1)\diam(I'(z_0))\leq
b_1(2\mu_{12}+1)d_{D'}(z'_0)\\
\nonumber &\leq&b_1(2\mu_{12}+1)(1+4c^2)^{\mu_1\mu_2}d_{D'}(x'),
\eeq
and thus \eqref{whuang-3-3} implies that
\beq\label{Thur-4}
\delta_{D'}(w',x')\leq 2 b_1(2\mu_{12}+1)(1+4c^2)^{\mu_1\mu_2} d_{D'}(x').
\eeq

Inequalities
\eqref{Thur-3} and \eqref{Thur-4} show that \eqref{monkey-1} is true by taking
$$b_2=2(1+4c^2)^{\mu_1\mu_2}\max\{(1+2c^2)^{\mu_1\mu_2}, b_1(2\mu_{12}+1)\}.$$

\subsubsection{\eqref{thm2-3} $\Longrightarrow$ \eqref{thm2-6}}\
\

For $x, w \in D$, we assume that $|x-w|\leq 8cd_D(x)$ and
$d_D(w)\leq 2cd_D(x)$. To prove the truth of the condition
\eqref{thm2-6} in the theorem, we consider two cases.

\setcounter{ca}{0}

\bca\label{Hmz-2}$|x-w|> \frac{1}{2}d_D(x)$.\eca

Let $y'_1\in D'$ be such that $|y'_1-w'|=\frac{1}{2}d_{D'}(w')$.
Then
$k_{D'}(y'_1,w')\leq
1,$ and so Theorem \Ref{ThmF} and (\ref{eq-2-9'})
imply that
$$\log\Big(1+\frac{\delta_D(y_1,w)}{d_{D}(w)}\Big)\leq k_{D}(y_1,w)
\leq \mu_2\max\Big\{k_D(y_1',w'), k_D(y_1',w')^{\frac{1}{\mu_2}}\Big\}\leq \mu_2,
$$
and so
\be\label{mlx-3-6-1}
\delta_D(y_1,w)\leq (\exp({\mu_2})-1)d_D(w).
\ee

It follows from the condition \eqref{thm2-3} in the theorem and
\cite{TV} that we may assume that $f$ is $\eta$-QS with
$\eta(t)=a\max\{t^{\frac{1}{\alpha}},t^{\alpha}\}$ for $t>0$, where the constants
$\alpha\in (0, 1]$ and $a$ depend only on $\varphi$ and the data
$$v=\Big(c, n, K,\frac{\diam(D)}{d_D(x_0)},
\frac{\diam(D')}{d_{D'}(f(x_0))} \Big).
$$
Then Lemma \ref{Sun-1} and the proved implication \eqref{thm2-4}
$\Longrightarrow$ \eqref{thm2-5} imply that the condition
\eqref{thm2-5} in the theorem holds. Hence we deduce from (\ref{mlx-3-6-1}) and
the assumption ``$d_D(w)\leq 2cd_D(x)$" that
\beq
\nonumber\frac{d_{D'}(w')}{2b_2d_{D'}(x')}&\leq&
\frac{\delta_{D'}(y'_1,w')}{\delta_{D'}(x',w')} \leq
\eta\Big(\frac{\delta_D(y_1,w)}{\delta_D(x,w)}\Big)\\
&\leq & \eta\Big(\frac{2(\exp({\mu_2})-1)d_D(w)}{d_D(x)}\Big)\\
\nonumber&<& a(2\exp({\mu_2}))^{\frac{1}{\alpha}}\max\bigg\{\Big(\frac{d_D(w)}{d_D(x)}\Big)^{\frac{1}{\alpha}},
\Big(\frac{d_D(w)}{d_D(x)}\Big)^{\alpha}\bigg\}\\
\nonumber&\leq&
a(4c\exp({\mu_2}))^{\frac{1}{\alpha}}\Big(\frac{d_D(w)}{d_D(x)}\Big)^{\alpha},
\eeq
which leads to
$$\frac{d_{D'}(w')}{d_D(w)}\leq
2ab_2(4c\exp({\mu_2}))^{\frac{1}{\alpha}}\frac{d_{D'}(x')}{d_D(x)}\cdot\Big(\frac{d_D(x)}{d_D(w)}\Big)^{1-\alpha},
$$
and by Theorem \Ref{mzxl-3-1}, we have
\beq\label{Sun-5}
a_f(w)\leq 2ab_2\mu_5^2(4c\exp({\mu_2}))^{\frac{1}{\alpha}}a_f(x)\Big(\frac{d_D(x)}{d_D(w)}\Big)^{1-\alpha}.
\eeq

\bca\label{Hmz-1}$|x-w|\leq \frac{1}{2}d_D(x)$.\eca

Obviously,
$$k_D(x,w)\leq \int_{[w,x]}\frac{|dz|}{d_D(z)}\leq
\frac{2|w-x|}{d_D(x)}\leq 3\log\frac{3d_D(w)}{d_D(x)},
$$
since $d_D(z)\geq d_D(x)-|x-z|\geq \frac{1}{2}d_D(x)$ for each
$z\in[w,x]$.
It follows from Theorem \Ref{ThmF} and (\ref{eq-2-9'}) that
\begin{eqnarray*}
\log\frac{d_{D'}(w')}{d_{D'}(x')}&\leq&
k_{D'}(x',w')\leq
\mu_2\max\Big\{k_D(x,w),(k_D(x,w))^{\frac{1}{\mu_2}}\Big\}\\ &\leq& 3\mu_2\log\frac{3d_D(w)}{d_D(x)},
\end{eqnarray*}
which implies that
$$\frac{d_{D'}(w')}{d_{D'}(x')}\leq \Big(\frac{3d_D(w)}{d_D(x)}\Big)^{3\mu_2},
$$
whence
$$\frac{d_{D'}(w')}{d_D(w)}\leq 3^{3\mu_2}\frac{d_{D'}(x')}{d_D(x)}\cdot \Big(\frac{d_D(x)}{d_D(w)}\Big)^{1-3\mu_2}.
$$
We infer from Theorem \Ref{mzxl-3-1} that
\beq\label{Sun-6}
a_f(w)\leq 2^{3\mu_2-\alpha}3^{3\mu_2}\mu_5^2a_f(x)\Big(\frac{d_D(x)}{d_D(w)}\Big)^{1-\alpha},
\eeq
since $d_D(w)\leq d_D(x)+|x-w|\leq \frac{3}{2}d_D(x)$.

We conclude from \eqref{Sun-5} and \eqref{Sun-6} that the condition
\eqref{thm2-6} in the theorem is true.

\begin{figure}[htbp]
\begin{center}

\input{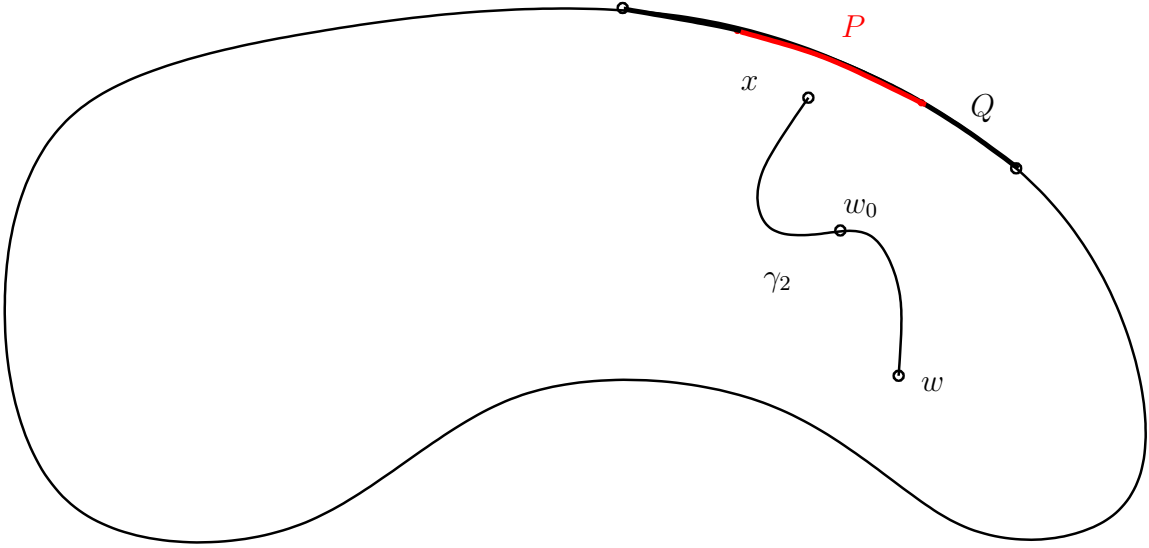}
\caption{The continua $P\subset Q$ and the points $x_1$, $x_2$,
$x_3$ and $x_4$.}

\label{hconj-new13}
\end{center}
\end{figure}

\subsubsection{ \eqref{thm2-6} $\Longrightarrow$ \eqref{thm2-7}}\
\

Suppose that the condition \eqref{thm2-6} in the theorem holds true
and that there are two points $x, w\in D$ such that $P\in \Phi(x)$, $Q\in\Phi(w)$ and $P\subset Q$
(see Figure \ref{hconj-new13}). To prove this implication, it suffices to show that
\beq\label{m-d-1}
\ds \frac{\diam(f(P))}{\diam(f(Q))}\leq b_4\Big(\frac{\diam(P)}{\diam(Q)}\Big)^{\alpha}.
\eeq

Apparently, $x=w$ implies $\Phi(x)=\Phi(w)$. In this case, $\eqref{m-d-1}$ is obvious. In the following,
we assume that $x\not= w$. We choose a
double $c$-cone arc $\gamma_2$ in $D$ joining $x$ and $w$ in $D$, and
use $w_0$ to denote the midpoint of $\gamma_2$ with respect to the
arc length (see Figure \ref{hconj-new13}).
Then the assumption $P\subset Q$ implies that
$$\mathbb{B}(x, 8cd_D(x))\cap \mathbb{B}(w, 8cd_D(w))\not=\emptyset,
$$
and so we have
\be\label{wangant-1}
\ell(\gamma_2)\leq c|x-w|\leq 8c^2(d_D(x)+d_D(w))
\ee
and
\be\label{wangant-2}
d_D(w_0)\geq \frac{1}{2c}\ell(\gamma_2)\geq \frac{1}{2c}|x-w|.
\ee

Since
$$d_D(x)\leq \diam(P)\leq \diam(Q)\leq 16cd_D(w),
$$
we get
\be\label{wangant-3}
|x-w|\leq 8c(d_D(x)+d_D(w))\leq 8c(1+16c)d_D(w).
\ee
Then it follows from Theorem \Ref{ThmF'}, (\ref{wangant-1}), (\ref{wangant-2}) and (\ref{wangant-3}) that
\begin{eqnarray*}
k_D(w,w_0)&\leq&\mu_1\log\Big(1+\frac{|w-w_0|}{\min\{d_D(w_0),d_D(w)\}}\Big)
\\ \nonumber&\leq& \mu_1\log\Big(1+\frac{\ell(\gamma_2)}{2\min\{d_D(w_0),d_D(w)\}}\Big)\\
\nonumber&\leq& \mu_1\log(1+4c^2+64c^3),
\end{eqnarray*}
which together with Theorem \Ref{ThmF} yield that
\begin{eqnarray*}
\log\frac{d_{D'}(w'_0)}{d_{D'}(w')}&\leq&
k_{D'}(w'_0,w')\\
&\leq &\mu_2\max\Big\{k_D(w_0,v),(k_D(w_0,v))^{\frac{1}{\mu_2}}\Big\}\\
&\leq& \mu_1\mu_2\log(1+4c^2+64c^3).
\end{eqnarray*}
Hence we have
\beq\label{samw-1}
\frac{d_{D'}(w'_0)}{d_{D'}(w')}\leq (1+4c^2+64c^3)^{\mu_1\mu_2}.
\eeq

Since
$$\max\{d_D(x),d_D(w)\}\leq 2cd_D(w_0)\;\;\mbox{and}\;\; |x-w_0|\leq cd_D(w_0),
$$
the condition \eqref{thm2-6} in the theorem and  Theorem \Ref{mzxl-3-1} lead to
$$\frac{d_{D'}(x')}{d_{D'}(w'_0)}\leq b_3\mu_5^2
\Big(\frac{d_D(x)}{d_D(w_0)}\Big)^{\alpha}\leq b_3(2c)^{\alpha}\mu_5^2
\Big(\frac{d_D(x)}{d_D(w)}\Big)^{\alpha},
$$
which, together with (\ref{samw-1}), shows that
$$\frac{d_{D'}(x')}{(1+4c^2+64c^3)^{\mu_1\mu_2}d_{D'}(w')}\leq \frac{d_{D'}(x')}{d_{D'}(w'_0)}
\leq b_3(32c^2)^{\alpha}\mu_5^2 \Big(\frac{\diam(P)}{\diam(Q)}\Big)^{\alpha},
$$
where in the last inequality, the following inequalities have been used: $d_D(x)\leq \diam(P)$ and $\diam(Q)\leq 16cd_D(w)$.
Hence by the condition \eqref{thm2-6} in the theorem, Lemmas \ref{lem5-A-11a} and \ref{lem5-A-11-1}, we have
\begin{eqnarray*}
\nonumber \frac{\diam(P')}{\mu_{12}\mu_{16}\diam(Q')}&\leq &
\frac{\diam(P')}{\mu_{16}d_{D'}(w')} \\
\nonumber &\leq &\frac{d_{D'}(x')}{d_{D'}(w')} \\
\nonumber &\leq & b_3(1+4c^2+64c^3)^{\mu_1\mu_2}(32c^2)^{\alpha}\mu_5^2\Big(\frac{\diam(P)}{\diam(Q)}\Big)^{\alpha},
\end{eqnarray*}
where in $\mu_{16}$, $\varepsilon_0=b_3$,
from which $\eqref{m-d-1}$ follows.

\subsubsection{\eqref{thm2-5} $\Longrightarrow$ \eqref{thm2-1}}\
\

To prove that $D'$ is a John domain, we let $\nu_1=\frac{\diam(D)}{d_D(x_0)}$. Then we have
\be\label{mxpl-5-2}
d_D(x_0)=\frac{1}{\nu_1}\diam(D).
\ee

For $z_1\in D\setminus\{x_0\}$, we let $\alpha_0$ be a double $c$-cone arc joining
$z_1$ and $x_0$ in $D$, and $y_0$ the midpoint of $\alpha_0$ with
respect to the arc length. Then we have

\setcounter{cl}{0}

\bcl\label{whuang5-1} For $v\in\alpha_0[y_0,x_0]$, we have
$d_D(v)\geq \frac{1}{2c}d_D(x_0)$.\ecl For a proof of this claim, we
consider two cases: $|x_0-v|\leq \frac{1}{2}d_D(x_0)$ and $|x_0-v|> \frac{1}{2}d_D(x_0)$.
If $|x_0-v|\leq \frac{1}{2}d_D(x_0)$ then
$$d_D(v)\geq d_D(x_0)-|x_0-v|\geq \frac{1}{2}d_D(x_0),
$$
and if $|x_0-v|>\frac{1}{2}d_D(x_0)$, then
$$d_D(v)\geq \frac{1}{c}\ell(\alpha_0[x_0, v]) \geq \frac{1}{c}|x_0-v|> \frac{1}{2c}d_D(x_0),
$$
since $\alpha_0$ is a double $c$-cone arc. The proof of Claim \ref{whuang5-1} is complete.\medskip

Further, we prove that $\alpha'_0$ is a carrot arc.

\bcl\label{whuang5-2-1}$\alpha'_0$ is a $b'_1$-carrot
$($diameter$)$ arc with the center $x'_0$, where
$$b'_1= 2\big(b_2\exp({c^2\nu_1\mu_2})+\exp({c^2\nu_1\mu_2})-1\big).
$$
\ecl

To prove this claim, it suffices to show that for $y\in\alpha_0$,
$$\diam(\alpha'_0[z'_1,y'])\leq b'_1d_{D'}(y').
$$
To this end, we separate the discussions into the following two cases.
On one hand, if $y\in \alpha_0[z_1,y_0]$, then there is some $y'_1$ in $\alpha'_0[z'_1,y']$
such that
$$|y'_1-y'|\geq \frac{1}{2}\diam(\alpha'_0[z'_1,y']).
$$
Since $|y_1-y|\leq\ell(\alpha_0[y_1,y])\leq cd_{D}(y)
$, we see from the condition \eqref{thm2-5} in the theorem
that
\be\label{whuang5-2} \nonumber
\diam(\alpha'_0[z'_1,y'])\leq 2|y'_1-y'|\leq
2\delta_{D'}(y'_1,y')\leq 2b_2d_{D'}(y')
\ee as required since $2b_2<b'_1$.

On the other hand, if $y\in \alpha_0[x_0,y_0]$, then for each $u\in\alpha_0[y_0,y]$, by
Claim \ref{whuang5-1} and (\ref{mxpl-5-2}), we have
$$k_D(y, u)\leq \int_{\alpha_0[y, u]}\frac{|dz|}{d_D(z)}\leq \frac{2c\ell(\alpha_0[y, u])}{d_D(x_0)}
\leq \frac{2c^2d_D(u)}{d_D(x_0)}\leq c^2\nu_1,
$$
since $\ell(\alpha_0[y, u])\leq cd_D(u)$. Hence Theorem \Ref{ThmF} and (\ref{eq-2-9'}) show that
\beq\label{king-kong-1}
\log\Big(1+\frac{|u'-y'|}{d_{D'}(y')}\Big)&\leq & k_{D'}(u',y')\leq
\mu_2\max\{k_D(u,y), (k_D(u,y))^{\frac{1}{\mu_2}}\}\\
\nonumber &\leq & c^2\nu_1\mu_2,
\eeq
whence
\be\label{whuang5-3}
|u'-y'|\leq \big(\exp({c^2\nu_1\mu_2})-1\big)d_{D'}(y').
\ee

In particular, \eqref{king-kong-1} also leads to
$$\log\frac{d_{D'}(y'_0)}{d_{D'}(y')}\leq k_{D'}(y'_0,y')\leq c^2\nu_1\mu_2,
$$
and so
\be\label{whuang5-4}
d_{D'}(y'_0)\leq \exp({c^2\nu_1\mu_2})d_{D'}(y').
\ee

By (\ref{whuang5-2}), (\ref{whuang5-3}) and (\ref{whuang5-4}), we
have
\beq\label{Sun-4}\nonumber
\diam(\alpha'_0[z'_1,y'])&\leq&
\diam(\alpha'_0[z'_1,y'_0])+\diam(\alpha'_0[y'_0,y'])\\
\nonumber&\leq&
2b_2d_{D'}(y'_0)+2\big(\exp({c^2\nu_1\mu_2})-1\big)d_{D'}(y')\\
\nonumber&\leq& 2\big(b_2\exp({c^2\nu_1\mu_2})+\exp({c^2\nu_1\mu_2})-1\big)d_{D'}(y')
\eeq as required. The proof
of Claim \ref{whuang5-2-1} is complete.
\medskip

The combination of Theorem \Ref{ThmF-1} and Claim \ref{whuang5-2-1}
shows that the condition \eqref{thm2-1} in the theorem is true. \qed

\subsection{
\eqref{thm2-1} $\Longrightarrow$ \eqref{thm2-8} $\Longrightarrow$
\eqref{thm2-9} $\Longrightarrow$ \eqref{thm2-10}  $\Longrightarrow$
\eqref{thm2-3}}\ \

The implications \eqref{thm2-1} $\Longrightarrow$ \eqref{thm2-8}
$\Longrightarrow$ \eqref{thm2-9} $\Longrightarrow$ \eqref{thm2-10}
follow from Theorems \Ref{Thm-1} with together with \cite[Lemmas 5.2 and 6.2]{JH}, and the implication
\eqref{thm2-10}  $\Longrightarrow$ \eqref{thm2-3} follows from
Theorem \ref{thm1}.\qed

\bigskip
\noindent {\bf Acknowledgements:}
The research of M. Huang and X. Wang was partly supported by NSFs of China (No. 11071063 and No. 11101138),
and was completed during the visit of these authors to IIT Madras.
They thank IIT Madras for the hospitality and the other support.

\end{document}